\documentclass[a4paper,11pt,reqno]{amsart}

\usepackage[ansinew]{inputenc}
\usepackage[english]{babel}
\usepackage{amsmath}
\usepackage{amssymb}
\usepackage{amsthm}
\usepackage{mathrsfs}
\usepackage[all]{xy}
\usepackage{verbatim}
\usepackage{arydshln}
\usepackage{multicol}
\usepackage{enumitem}
\usepackage{bm}
\usepackage{bbm}

\usepackage[bookmarks=true,bookmarksopen=false]{hyperref}
\hypersetup{
   pdftitle = {},
   pdfauthor = {Matthew Gelvin},
   pdfpagemode = {UseOutlines},
   pdfstartview = {FitH},
   pdfborder = {0 0 0},
   backref = {true},
   colorlinks = {true},
   urlcolor = {blue},
   citecolor = {purple},
   linkcolor = {blue},
   pdftoolbar = {true},
}

\usepackage{mathtools}
\usepackage{ifthen}
\usepackage{tikz}
\usepackage{todonotes}
\usepackage{bookmark}

\allowdisplaybreaks

\theoremstyle{plain}
\newtheorem{theorem}{Theorem}
\newtheorem{prop}[theorem]{Proposition}
\newtheorem{lemma}[theorem]{Lemma}
\newtheorem{cor}[theorem]{Corollary}
\newtheorem{conj}{Conjecture}

\theoremstyle{definition}
\newtheorem{definition}[theorem]{Definition}

\theoremstyle{remark}

\DeclareMathOperator{\Aut}{Aut}

\DeclareMathOperator{\End}{End}
\DeclareMathOperator{\Stab}{Stab}

\DeclareMathOperator{\Hom}{Hom}
\DeclareMathOperator{\Iso}{Iso}

\DeclareMathOperator{\Syl}{Syl}

\DeclareMathOperator{\tra}{tr}
\DeclareMathOperator{\br}{br}

\newcommand{\id}{\mathrm{id}}

\renewcommand{\tilde}{\widetilde}

\def\cA{\mathcal A}\def\cB{\mathcal B}
\def\cF{\mathcal F}\def\cG{\mathcal G}
\def\cL{\mathcal L}
\def\cO{\mathcal O}\def\cP{\mathcal P}
\def\cS{\mathcal S}

\def\fL{\mathfrak L}

\def\td{\mathtt{d}}\def\tg{\mathtt{g}}\def\th{\mathtt{h}}
\def\tp{\mathtt{p}}\def\tq{\mathtt{q}}
\def\ts{\mathtt{s}}
\def\tD{\mathtt{D}}\def\tG{\mathtt{G}}\def\tH{\mathtt{H}}
\def\tP{\mathtt{P}}\def\tQ{\mathtt{Q}}
\def\tR{\mathtt{R}}\def\tS{\mathtt{S}}\def\tT{\mathtt{T}}

\def\fF{\mathfrak F}

\def\kk{\mathbbm{k}}

\def\id{\mathrm{id}}

\usetikzlibrary{arrows,matrix,decorations,positioning,calc}
\tikzset{dot/.style={circle,fill=black,thick,inner sep=0pt,minimum size=1mm,draw}}
\tikzset{arrow/.style={semithick,>=stealth',shorten >=1pt,shorten <=1pt}}
\tikzset{equal/.style={arrow,double distance=2pt}}

\title[Source algebras and characteristic bisets]{Variants on a conjecture relating block source algebras to characteristic bisets}
\author{Laurence Barker and Matthew Gelvin}
\date{\today}

\begin{document}

\maketitle
\begin{abstract}
Given a block of a finite group, any source algebra has a basis invariant under the multiplicative actions of the defect group.  Is such a basis a characteristic biset of the block fusion system?  If the basis can be chosen to consist entirely of units, the question is answered in the affirmative.  We prove the equivalence of several reformulations of this stronger condition on the source algebra.
\end{abstract}

\section{Introduction}

Many important invariants of a block of a finite group, perhaps all invariants that ought to be deemed local, are determined by its source algebra.  This algebra is interior for the defect group, with the special property that it possesses a linear basis permuted by the left and right group actions.  We consider the conjecture that any such basis is a characteristic biset of the block fusion system.  We also examine the stronger condition that an invariant basis can be chosen to consist entirely of units.  

Let $\tG$ be a finite group and $\cO$ a complete local noetherian domain whose residue field $\kk$ is algebraically closed  of prime characteristic $p$, allowing for the possibility that $\cO=\kk$.  In the group algebra $\cG:=\cO\tG$, each block idempotent $b$ has a   defect group $\tD$,  unique up to $\tG$-conjugacy. The block algebra $\cB:=b.\cG$ is an interior $\tD$-algebra, so we may view it as an $(\cO\tD,\cO\tD)$-bimodule. As such, $\cB$ has an $\cO$-basis $Y$ that is invariant under the natural $(\tD,\tD)$-biaction:  For all $\td\in \tD$, we have $\td\cdot Y=Y=Y\cdot \td$.  This makes $Y$  a $(\tD,\tD)$-biset, and a different choice of invariant basis of $\cB$ yields an isomorphic biset.  

The block idempotent $b$ also gives rise to the block fusion system $\cF_\tD(b)$ on $\tD$,  uniquely determined up to isomorphism of fusion systems.  We begin by asking what relationship there may be between this fusion system and the $(\tD,\tD)$-biset $Y$.  

In their unpublished notes, Linckelmann and Webb introduced ``characteristic bisets'' as a bridge between the worlds of fusion  and bisets.  If $\cF$ is a fusion system on the $p$-group $\tS$, an $\cF$-characteristic biset is an $(\tS,\tS)$-biset that determines and (to some extent, cf. \cite{GelvinReehMinimalCharacteristicBisets,GelvinReehYalcin}) is determined by $\cF$, according to certain properties enumerated in Definition \ref{def:F_characteristic}.  Three of the defining properties are worth naming here:  A characteristic biset must be ``$\cF$-generated'' (the `conjugation' maps induced by point-stabilizers must be morphisms of $\cF$), ``$\cF$-stable'' (so that the biset is unchanged, up to isomorphism, when twisted by a morphism of $\cF$), and it must satisfy a Sylow $p$-primality condition. Roughly speaking, an $\cF$-characteristic biset is a stand-in for an ambient finite group that contains $\tS$ as a Sylow $p$-subgroup and induces $\cF$. Even when no such finite group exists, it is still the case that every saturated fusion system has a characteristic biset \cite[Proposition 5.5]{BLO2}.  Conversely, the existence of a characteristic biset implies the saturation of the fusion system (cf. \cite[Proposition 21.9]{PuigBook} and \cite[\S 6]{RagnarssonStancuIdempotents}).

 Is $Y$ an $\cF_\tD(b)$-characteristic biset?  One quickly sees that it cannot be so in general.  If $\tS$ is a Sylow $p$-subgroup of $\tG$, the well-known formula (cf. \cite[Theorem 1]{BrauerNotesOnRepresentations}) for the $p$-part of the $\cO$-rank of $\cB$ implies that $(|Y|/|\tD|)_p=(|\tG|_p/|\tD|)^2$. This expression plays the role of the `index' of $\tD$ in $Y$, so we roughly have that $\tD$ is not a `Sylow $p$-subgroup of $Y$' unless $\tD$ is a Sylow $p$-subgroup of $\tG$.  Thus the Sylow $p$-primality condition may fail for an invariant basis of a block algebra.
 
There is another, perhaps more fundamental, failure of $Y$ to be $\cF_\tD(b)$-characteristic.  It is not hard to see that for $Y$ an $\cO$-invariant basis of $\cB$, the induced `conjugation by an element of $Y$' maps (i.e., an assignment of the form $\tp\mapsto \tp'$ for $\tp\in\tP\leq\tS$ and $\tp'\in\tS$ such that there exists $y\in Y$ with $y\cdot\tp=\tp'\cdot y$) are the same as those maps induced by conjugation in $\tG$.  In other words, $Y$ can only be $\cF_\tD(b)$-generated if $\cF_\tD(b)$ a full subcategory of the  fusion system induced by $\tG$, which is rarely the case.  Thus an invariant basis for a block algebra need not be $\cF_\tD(b)$-generated.

The heart of the issue is that, while a block does determine a fusion system  up to isomorphism, there are certain choices that must be made to specify it uniquely at the level of equality.  We have already implicitly settled on a representative defect group $\tD$ from its $\tG$-conjugacy class; this names the  objects of $\cF_\tD(b)$, but not the morphisms.  To fully specify the block fusion system, we must also choose a ``maximal Brauer pair'' associated to $b$.  Equivalently, we must choose a primitive idempotent $\ell$ in the $\tD$-fixed subalgebra $\cB^\tD$ that is ``local'' in the sense that $\br_\tD(\ell)\neq 0$ (cf. \S\ref{subsec:twisted_Brauer_quotients}). We will continue to write $\cF_\tD(b)$ for what might more properly be denoted $\cF_\tD(\ell)$, with the understanding that a particular choice of block fusion system has been made.  See \S\ref{sec:Block_fusion_and_local_categories} for a review of this material.

The idempotent $\ell$ is known as a \emph{source idempotent} of $b$, and it gives rise to the \emph{source algebra} $\cS:=\ell.\cB.\ell$, which is again an interior $\tD$-algebra with $(\tD,\tD)$-invariant $\cO$-basis $X$. As $\cS$ is a corner algebra of $\cB$, it follows that we may view $X$ as a $(\tD,\tD)$-subbiset of $Y$.  The failures of $Y$'s being $\cF_\tD(b)$-characteristic stem from there being either too many $(\tD,\tD)$-orbits or $(\tD,\tD)$-orbits of the wrong form, so we might hope that moving from $Y$ to $X$ will eliminate just the right orbits to obtain a $\cF_\tD(b)$-generated biset that satisfies the Sylow $p$-primality condition.  This turns out to be the case, but the price we pay is that we may lose $\cF_\tD(b)$-stability.  Nevertheless, we are unaware of any general reason that $\cF_\tD(b)$-stability should fail, and we have been unable to find an example where it does.  We are thus led to the motivating conjecture:

\begin{conj}\label{conj:main}
A $(\tD,\tD)$-invariant $\cO$-basis $X$ of a source algebra $\cS$ is $\cF_\tD(b)$-characteristic.
\end{conj}

We do not prove Conjecture \ref{conj:main} in this paper.  Instead, our aim is to describe certain structural properties of $\cS$ that will imply this result, should they be verified.  We identify three such sufficient criteria:  the existence of a $(\tD,\tD)$-invariant $\cO$-basis of $\cS$ that is contained in $\cS^\times$ (termed ``$\cS$ has a unital $(\tD,\tD)$-invariant $\cO$-basis''), the coherent regularity of certain ``twisted'' Brauer quotients  of $\cS$ as modules over the regular Brauer quotient algebras (``$\cS$ has all twisted units''), and the high degree of symmetry of the embedding  $\cS\subseteq\cB$ (``$\cS$ is balanced in $\cB$'').  The first two of these explicitly relate the multiplication of $\cS$ to properties of $X$ that would guarantee $\cF_\tD(b)$-stability, while the third examines the multiplicities of certain idempotents to ask when the $\cF_\tD(b)$-stability of $\cB$  is  inherited by $\cS$.  We outline these criteria now.

When first considering the question of the $\cF_\tD(b)$-stability of $X$, we were reminded of the work of Ragnarsson and Stancu \cite{RagnarssonStancuIdempotents}, which gave a necessary and sufficient condition on an $(\tS,\tS)$-biset $\Omega$ to guarantee it be a characteristic idempotent for \emph{some} fusion system.  The ``Frobenius reciprocity'' property they identify asks that certain  $\tS\times \tS\times \tS$-structures one could put on $\Omega\times \Omega$ be isomorphic. As $X\times X$ can be identified with a $\cO$-basis for $\cS\otimes\cS$ via the map $(x,x')\mapsto x\otimes x'$, it was natural to consider the corresponding $\cO(\tD\times \tD\times \tD)$-module structures one could place on $\cS\otimes\cS$, and ask whether they could be seen to yield isomorphic modules.  The map $X\times X\to\cS\otimes\cS:(x,x')\mapsto x\otimes x.x'$ has the right equivariance properties, but its image is only guaranteed to be an $\cO$-basis if the elements of $X$ are all units in $\cS$.  If such a unital basis could be found, then $X$ would be an $\cF$-characteristic biset for some saturated fusion system $\cF$ on $\tD$, and the established properties of $X$ would imply $\cF=\cF_\tD(b)$.

\begin{theorem}\label{thm:unital_bases_imply_characteristic}
If the source algebra $\cS$ has a $(\tD,\tD)$-invariant $\cO$-basis $X$ such that $X\subseteq\cS^\times$, then $X$ is  $\cF_\tD(b)$-characteristic.
\end{theorem}

There are two points worth noting.  The first regards authorship: This result was known to Linckelmann and Webb in their unpublished work on characteristic bisets, and it later appeared as \cite[Proposition 8.7.11]{LinckelmannBookII}.  We justify our inclusion of Theorem \ref{thm:unital_bases_imply_characteristic} here by noting that we actually prove a slight generalization (Lemma \ref{lem:units_in_fixed_point_submodules} and Proposition \ref{prop:basis_of_units_implies_conj1}), which will be key to proving Theorems \ref{thm:twisted_units_imply_main} and \ref{thm:balance_implies_main_conjecture}.

The second point is that, after realizing the importance of the existence of a unital basis, it was easy to find a more direct proof that did not rely on the Frobenius reciprocity condition.  This is the argument  presented in \S\ref{sec:unital_bases}.  We mention the history  out of the sense that there is untapped potential in the notions introduced by Ragnarsson and Stancu, and offer the above as  evidence of this claim.

In any event, the bold may now be led to

\begin{conj}\label{conj:main_strong_i}
The source algebra $\cS$ has a unital $(\tD,\tD)$-invariant $\cO$-basis.
\end{conj}

Our second approach to Conjecture \ref{conj:main} is more direct, using the formulation of $\cF_\tD(b)$-stability in terms of the orders of certain fixed-point subsets of $X$.  Explicitly,  for any $\tP\leq\tD$ and group injection $\varphi:\tP\hookrightarrow\tD$, one may form the ``twisted diagonal subgroup'' $\Delta(\varphi,\tP):=\{(\varphi(\tp),\tp)\ \big|\ \tp\in \tP\}\leq\tD\times\tD$. Then  $X$ is $\cF_\tD(b)$-stable if and only if for all $\varphi\in\Iso_{\cF_\tD(b)}(\tP,\tQ)$, the number points of $X$ fixed by $\Delta(\varphi,\tP)$, $\Delta(\tP)$, and $\Delta(\tQ)$ are all equal.  

By viewing the source algebra $\cS$ as a $k(\tD\times \tD)$-module we may speak of the Brauer quotient $\cS(\tR)$ for any subgroup $\tR\leq \tD\times \tD$, and the existence of a $(\tD,\tD)$-stable $\cO$-basis allows us to equate the $\mathbbm{k}$-dimension of $\cS(\tR)$ with the number of $\tR$-fixed-points of $X$.  In particular, the $\cF_\tD(b)$-stability of $X$ is equivalent to the equality 
\[
\dim_\kk\cS(\tQ)=\dim_\kk\cS(\Delta(\varphi,\tP))=\dim_\kk\cS(\tP)
\]
for all isomorphisms $\varphi\in\Iso_{\cF_\tD(b)}(\tP,\tQ)$.  Here $\cS(\tP)$ and $\cS(\tQ)$ are the standard Brauer quotients obtained by viewing $\cS$ as a $\tD$-algebra, and $\cS(\Delta(\varphi,P))$ will be known as the ``twisted Brauer quotient'' of $\cS$ at $\varphi$, as described in \S\ref{subsec:twisted_Brauer_quotients}.  

The algebra structure of $\cS$ descends to $\cS(\tQ)$ and $\cS(\tP)$, making each  an algebra in its own right.  Moreover, $\cS$-multiplication induces an $(\cS(\tQ),\cS(\tP))$-bimodule structure on  $\cS(\Delta(\varphi,\tP))$.  The most natural way to obtain our desired equality of Brauer quotient dimensions would be to show that $\cS(\Delta(\varphi,\tP))$ is  regular as both a left $\cS(\tQ)$- and right $\cS(\tP)$-module.  Even better would be to show that a single element of $\cS(\Delta(\varphi,\tP))$ is a  generator for both regular module structures, and best of all would be if the set of all such generators were compatible with the composition of  isomorphisms  in $\cF_\tD(b)$.  We call such a simultaneous generator a ``twisted unit'' of $\cS$ at $\varphi$, and if a coherent collection of twisted units exists we say that $\cS$   ``has all twisted units.''  See \S\ref{sec:twisted_units}.

\begin{theorem}\label{thm:twisted_units_imply_main}
If the source algebra $\cS$ has all twisted units, then any $(\tD,\tD)$-invariant $\cO$-basis of $\cS$ is $\cF_\tD(b)$-charactersitic.
\end{theorem}

\begin{conj}\label{conj:main_strong_ii}
The source algebra $\cS$ has all twisted units.
\end{conj}

Our final approach makes note of the $\cF_\tD(b)$-stability of the block algebra $\cB$ and asks when this property is inherited by the source algebra $\cS=\ell.\cB.\ell$. In \S\ref{sec:balanced_algebras} we identify a sufficient a condition, which we express by saying $\cS$ is ``balanced'' in $\cB$.  Roughly speaking, balance asserts that for every $\varphi\in\Iso_{\cF_\tD(b)}(\tP,\tQ)$, there should be some unit in the $\Delta(\varphi,\tP)$-fixed submodule of $\cB$ that induces a bijection between the local points of $\tP$ and $\tQ$ on $\cS$, and  that moreover this bijection should preserve the multiplicities of local points.
We also give an intrinsic description of the balance condition without reference to the ambient block algebra.  We prove:

\begin{theorem}\label{thm:balance_implies_main_conjecture}
If the source algebra $\cS$ is balanced in the block algebra $\cB$, then any $(\tD,\tD)$-invariant $\cO$-basis of $\cS$ is $\cF_\tD(b)$-characteristic.
\end{theorem}

\begin{conj}\label{conj:main_strong_iii}
The source algebra $\cS$ is balanced in the block algebra $\cB$.
\end{conj}

At this point we have  Conjecture \ref{conj:main} and the three stronger Conjectures \ref{conj:main_strong_i}, \ref{conj:main_strong_ii},  and \ref{conj:main_strong_iii}.  It is not hard to see that in all cases where our motivating conjecture is known to hold each of the stronger ones does as well, but they appear otherwise unrelated at first glance.  Our main result is that the hypotheses of the stronger conjectures are in fact equivalent.

\begin{theorem}\label{thm:main}
Let $\cS$ be a source algebra of the block algebra $\cB$ with defect group $\tD$.  The following are equivalent:
\begin{enumerate}
\item\label{thm:main_i} $\cS$ possesses a unital $(\tD,\tD)$-invariant $\cO$-basis.
\item\label{thm:main_ii} $\cS$ has all twisted units.
\item\label{thm:main_iii} $\cS$ is balanced in $\cB$.
\end{enumerate}
\end{theorem}

We actually prove more than is stated in Theorem \ref{thm:main}, as this equivalence holds for all sufficiently structured $\cO$-algebras, which we term ``divisible'' (see \S\ref{subsec:twisted_Brauer_quotients}).  We find it useful to work in this greater level of generality, as  a companion paper will consider and confirm analogues of Conjectures \ref{conj:main} through \ref{conj:main_strong_iii}  for certain ``almost-source algebras'' (cf. \cite{LinckelmannTrivialSource}) of  $p$-solvable groups.  (Indeed, our conjectures are posed in the most ambitious manner conceivable, and Linckelmann has suggested that perhaps the proper statements should replace ``any source algebra of a block satisfies\ldots'' with ``any block has an almost-source algebra that satisfies\ldots.'') Almost-source algebras are almost source algebras in terms of the structural properties of interest to us, so our more general viewpoint is taken with this application in mind.

For this reason, most statements of this paper will be made in terms of a finite dimensional $\cO$-free $\cO$-algebra $\cA$, which will be acted on by a finite $p$-group $\tS$.  We reserve the symbols ``$\cB$,'' ``$\cS$,'' and ``$\tD$'' for the block algebra, source algebra, and their common defect group, respectively. We will generally only make explicit references to the latter when we need to verify that certain desired properties are satisfied by our motivating example.

The proof of Theorem \ref{thm:main} is split into multiple parts throughout this paper.  In the case of divisible algebras, the equivalence \ref{thm:main_i}$\Leftrightarrow$\ref{thm:main_ii}  appears as Corollary \ref{cor:1_2_equiv}, while \ref{thm:main_i}$\Leftrightarrow$\ref{thm:main_iii} is Theorem \ref{thm:balance_and_unital_basis}.  The connection to our conjectures comes from Theorem \ref{thm:source_algebras_are_divisible}, where we prove that source algebras are divisible, and so these more general results apply.  

This paper is organized as follows:
\begin{enumerate}
\item[\S\ref{sec:preliminaries}] reminds the reader of notions from the literature and establishes our notation.
\begin{enumerate}
\item[\S\ref{subsec:biset_and_fusion_intro}] defines $\cF$-characteristic bisets and related notions.
\item[\S\ref{subsec:permutation_algebra_biset_bases}] summarizes the foundational work of Green \cite{GreenBlocksOfModularRepresentations} on permutation $\tS$-algebras.  In particular, we note that if the $\tS$-action on $\cA$ is interior, then all $(\tS,\tS)$-invariant $\cO$-bases of $\cA$ determine the same $(\tS,\tS)$-biset, up to isomorphism.
\item[\S\ref{subsec:twisted_Brauer_quotients}] recalls the Brauer quotient construction and specializes to the twisted Brauer quotients that serve as a bridge between the module structure of $\cA$ and the $(\tS,\tS)$-invariant $\cO$-basis $Y$.  This section also introduces the notion of a ``divisible'' $\tS$-algebra, where the point-stabilizers of $Y$ form a fusion system on $\tS$, and observes that the existence of a unital $\cO$-basis implies divisibility.
\item[\S\ref{subsec:state_of_the_art}] is an overview of the current status of Conjecture \ref{conj:main}.
\end{enumerate}
\item[\S\ref{sec:Block_fusion_and_local_categories}] reviews the definition of the block fusion system $\cF_\tD(b)$ and Puig's theory of  local categories.  The latter is of particular interest, as it relates more directly to the $(\tD,\tD)$-invariant $\cO$-basis of $\cS$ than the block fusion system itself. Here we also record a number of technical lemmas that will be needed in \S\S\ref{sec:unital_bases}--\ref{sec:balanced_algebras}, and prove that source algebras of blocks are divisible.
\item[\S\ref{sec:unital_bases}] proves that Conjecture \ref{conj:main_strong_i} implies Conjecture \ref{conj:main}.  Here we also discuss a necessary and sufficient condition for the existence of a unital $(\tS,\tS)$-invariant $\cO$-basis of $\cA$, stated in terms of the fixed-point subsets of the unit group $\cA^\times$.
\item[\S\ref{sec:twisted_units}] introduces the notion of twisted units. We prove that $\cA$ has a unital $(\tS,\tS)$-invariant $\cO$-basis if and only if it has all twisted units.
\item[\S\ref{sec:balanced_algebras}] defines the balance condition, both as a relationship between an $\tS$-algebra $\widehat{\cA}$ and a corner subalgebra $\cA=\ell.\widehat{\cA}.\ell$, and as an intrinsic proptery of $\cA$ itself. We prove that the existence of a unital $(\tS,\tS)$-invariant $\cO$-basis of $\cA$ is equivalent to intrinsic balance.  This completes the proof of Theorem \ref{thm:main}, and the paper.
\end{enumerate}

\numberwithin{theorem}{section}

\section{Preliminaries}\label{sec:preliminaries}

\subsection{Bisets and fusion systems}\label{subsec:biset_and_fusion_intro}  Let $\tG$ and $\tH$ be finite groups.

A \emph{$(\tG,\tH)$-biset} is a finite set $\Omega$ that is simultaneously a left $\tG$-set and a right $\tH$-set, where each action respects the other:
\[
\tg\cdot(\omega\cdot \th)=(\tg\cdot\omega)\cdot \th\qquad\textrm{for all }\tg\in \tG,\ \th\in \tH,\textrm{ and }\omega\in\Omega.
\]
The \emph{stabilizer} of a point $\omega\in\Omega$ is
\[
\Stab_{(\tG,\tH)}(\omega)=\Stab(\omega):=\left\{(\tg,\th)\in \tG\times \tH\ \big|\ \tg\cdot\omega=\omega\cdot \th\right\},
\]
which is clearly a subgroup of $\tG\times \tH$.  Dually, for $\mathtt{A}\leq \tG\times \tH$, the \emph{$\mathtt{A}$-fixed-points} of $\Omega$ are
\[
\Omega^\mathtt{A}:=\left\{\omega\in\Omega\ \big|\ \mathtt{a}_\tG\cdot\omega=\omega\cdot \mathtt{a}_\tH\ \forall\ (\mathtt{a}_\tG,\mathtt{a}_\tH)\in \mathtt{A}\right\}.
\]

$\Omega$ is \emph{bifree} if the left and right group actions are individually free.  Goursat's Lemma implies that each point-stabilizer in a bifree biset is of the form
\[
\Delta(\varsigma,\mathtt{K}):=\left\{(\varsigma(\mathtt{k}),\mathtt{k})\ \big|\ \mathtt{k}\in \mathtt{K}\right\}
\]
for some subgroup $\mathtt{K}\leq \mathtt{H}$ and group monomorphism $\varsigma:\mathtt{K}\hookrightarrow \mathtt{G}$.  Groups of the form $\Delta(\varsigma,\mathtt{K})$ will be called \emph{twisted diagonal subgroups} of $\tG\times\tH$.

Observe that if $\Omega$ is bifree and $\mathtt{A}\leq \tG\times \tH$, then $\Omega^\mathtt{A}\neq\emptyset$ only if $\mathtt{A}$ is a twisted diagonal subgroup.  We  adopt the notation ${^\varsigma\Omega ^\mathtt{K}}$ for $\Omega^{\Delta(\varsigma,\mathtt{K})}$ in this case.  Moreover, when $\tG=\tH$, $\mathtt{K}\leq \tG$, and $\iota_\mathtt{K}^\tG:\mathtt{K}\hookrightarrow \tG$ is the natural inclusion, we will write $\Omega^\mathtt{K}$ for ${^{\iota_\mathtt{K}^\tG}\Omega^\mathtt{K}}=\Omega^{\Delta(\iota_\mathtt{K}^\tG,\mathtt{K})}$.

The \emph{opposite biset} of the $(\tG,\tH)$-biset $\Omega$ is the $(\tH,\tG)$-biset $\Omega^\circ$ with the same underlying set and action given by
\[
\th\odot\omega\odot \tg:=\tg^{-1}\cdot\omega\cdot \th^{-1}.
\]
If $\tG=\tH$ and $\Omega\cong\Omega^\circ$ as $(\tG,\tG)$-bisets, we say that $\Omega$ is \emph{symmetric}.

Bisets arise naturally in the study of fusion systems.  Recall that a \emph{fusion system} on the  $p$-group $\tS$ is a category $\cF$ whose objects are the subgroups of $\tS$ and whose morphisms are  certain collections of injective group maps.  A \emph{realizable} fusion system is one of the form $\cF_\tS(\tG)$, where $\tG$ is a finite group, $\tS\in\Syl_p(\tG)$, and the morphisms are  $\tG$-conjugations.  The block fusion system of interest to us, $\cF_\tD(b)$, will be recalled in greater detail in \S\ref{sec:Block_fusion_and_local_categories}.

Brauer's Third Main Theorem tells us that every realizable fusion system is the fusion system of a block.  Both of these are examples of \emph{saturated} fusion systems, which satisfy additional axioms meant to mimic the realizable case.  We have no need of the precise definition of saturation in this paper, and will only use a consequence of a block fusion system's saturation in the proof of Theorem \ref{thm:source_algebras_are_divisible}.  For the general theory, the reader may refer to Puig's work \cite{PuigFrobeniusCategories} or the more recent \cite{AKO}; the proof that block fusion systems are saturated can be found in Part IV of the latter.

The connection between fusion and bisets comes from considering realizability from a different perspective:  Instead of asking how $\tG$ acts on the subgroups of $\tS$ by conjugation, we consider how $\tS$ acts on $\tG$ by left and right multiplication.  The resulting $(\tS,\tS)$-biset $_\tS\tG_\tS$ is not as structurally rich as the group $\tG$,  yet it still contains enough information to determine the fusion system $\cF_\tS(\tG)$. We now recall how this is achieved.

Fix a fusion system $\cF$ on $\tS$, and let $\Omega$ be a symmetric bifree $(\tS,\tS)$-biset.   We say that $\Omega$ is \emph{$\cF$-generated} if every $\omega\in\Omega$ has point-stabilizer of the form $\Delta(\varphi,\tP)$ with $\varphi\in\Hom_\cF(\tP,\tS)$.  Conversely, $\Omega$ is \emph{$\cF$-stable} if for every $\varphi\in\Hom_\cF(\tP,\tS)$, there is an isomorphism of $(\tP,\tS)$-bisets ${_\tP\Omega_\tS}\cong{_\tP^\varphi\Omega_\tS}$, where the $(\tP,\tS)$-biset structure of $_\tP\Omega_\tS$ obtained by restriction and that of $^\varphi_\tP\Omega_\tS$ by restriction and twisting the left $\tP$-action along $\varphi$:  For  $\tp\in \tP$, $\ts\in \tS$, and $\omega\in\Omega$,  set $\tp\odot\omega\odot \ts:=\varphi(\tp)\cdot\omega\cdot \ts$.  Observe that $\cF$-stability is equivalent to requiring the equality of fixed-point orders $\left|{^\varphi\Omega^\tP}\right|=\left|\Omega^\tP\right|$ for all $\tP\leq \tS$ and $\varphi\in\Hom_\cF(\tP,\tS)$.

\begin{definition}\label{def:F_characteristic}
The $(\tS,\tS)$-biset $\Omega$ is \emph{$\cF$-characteristic} if 
\begin{enumerate}
\item\label{def:F_characteristic_i} $\Omega$ is bifree,
\item\label{def:F_characteristic_ii} $\Omega$ is symmetric,
\item\label{def:F_characteristic_iii} $\Omega$ is $\cF$-generated,
\item\label{def:F_characteristic_iv} $\Omega$ is $\cF$-stable, and
\item\label{def:F_characteristic_v} $|\Omega|/|\tS|$ is prime to $p$.
\end{enumerate}
\end{definition}

Thus proving Conjecture \ref{conj:main} is just a checklist verification of the conditions of Definition \ref{def:F_characteristic} for a $(\tD,\tD)$-invariant $\cO$-basis of $\cS$ and block fusion system $\cF_\tD(b)$, cf. \S\ref{subsec:state_of_the_art}.  

\subsection{Interior permutation algebras and biset bases}\label{subsec:permutation_algebra_biset_bases}

Let $\tG$ be a finite group and $\cA$ an $\cO$-algebra that is free as an $\cO$-module.

A \emph{$\tG$-algebra} structure on $\cA$ is an action of $\tG$ on $\cA$ by algebra automorphisms.   If $\cA$ possesses a $\cO$-basis that is invariant under the $\tG$-action, we say that $\cA$ is a \emph{permutation} $\tG$-algebra.  More generally, $\cA$ is a \emph{$p$-permutation} $G$-algebra if restricting the action of $\tG$ to any $\tS\in\Syl_p(\tG)$ yields a permutation $\tS$-algebra.

Suppose now that we are given a group map $\tG\to\cA^\times$, allowing us to view a homomorphic image $\overline{\tG}$ of $\tG$ as a subgroup of the unit group of $\cA$.  In this case we say that $\cA$ is an \emph{interior $\tG$-algebra}, and we may view $\cA$ as a $(\cO\tG,\cO\tG)$-bimodule with action induced by $\tg_1\cdot a\cdot \tg_2=\overline{\tg_1}.a.\overline{\tg_2}$.  In particular, we may view an interior $\tG$-algebra as a $\tG$-algebra via conjugation:  For $\tg\in \tG$ and $a\in\cA$, the action of $\tg$ on $a$ is given by $^\tg a:=\tg\cdot a\cdot \tg^{-1}$.  

The interior $\tG$-algebra $\cA$ is a \emph{bipermutation  $\tG$-algebra} if it possesses an $\cO$-basis that is invariant under the left and right actions of $\tG$.  Such a basis will be called \emph{$(\tG,\tG)$-invariant}, and is naturally a $(\tG,\tG)$-biset.  Moreover, we say that $\cA$ is a \emph{bifree $\tG$-algebra} if it is a bipermutation $\tG$-algebra and any of its $(\tG,\tG)$-invariant $\cO$-bases are bifree as $(\tG,\tG)$-bisets.  As above, we define $\cA$ to be a \emph{$p$-bipermutation $\tG$-algebra} if restriction of the group map $\tG\to\cA^\times$ to any $\tS\in\Syl_p(\tG)$ yields a bipermutation $\tS$-algebra.

The notions of \emph{$\tG$-module}, \emph{permutation $\tG$-module}, and \emph{$p$-permutation $\tG$-module} are defined in the obvious way.  Our main object of study, the source algebra $\cS$ of a $p$-block algebra $\cB$ with defect group $\tD$, is a permutation interior $\tD$-algebra. As such, the underlying $\cO$-module structure of $\cS$ is that of a permutation $\cO(\tD\times \tD)$-module.  As $\tD\times \tD$ is a $p$-group, we may make use of Green's Indecomposability Criterion (cf. \cite{GreenBlocksOfModularRepresentations}).  The key results for us are:
\begin{prop} 
Let $\tS$ be a $p$-group with subgroup $\tP\leq \tS$, $[\tS/\tP]$ the corresponding transitive $\tS$-set, and $\cO[\tS/\tP]$ the permutation $\tS$-module with $\cO$-basis $[\tS/\tP]$.
\begin{enumerate}
\item $\cO[\tS/\tP]$ is indecomposable as an $\tS$-module.
\item If $\cO[\tS/\tP]\cong \cO[\tS/\tQ]$ as $\tS$-modules, then $[\tS/\tP]\cong[\tS/\tQ]$ as $\tS$-sets.
\end{enumerate}
\end{prop}
\begin{cor}\label{cor:basis_well_defined}
Let $M$ be a permutation $\tS$-module with $\tS$-invariant $\cO$-basis $Y$.
\begin{enumerate}
\item\label{cor:basis_well_defined_i}  If $Y'$ is another $\tS$-invariant $\cO$-basis of $M$, then $Y\cong Y'$ as $\tS$-sets.
\item\label{cor:basis_well_defined_ii} If $N$ is a direct summand of $M$, then $N$ is also a permutation $\tS$-module.
\item\label{cor:basis_well_defined_iii} If the direct summand $N$ of $M$ has $\tS$-invariant $\cO$-basis $Z$, then $Z$ is isomorphic to an $\tS$-subset of $Y$.
\end{enumerate}
\end{cor}

In particular, the $(\tD,\tD)$-invariant $\cO$-basis $X$ of $\cS$  is well-defined up to isomorphism of $(\tD,\tD)$-bisets, and so Conjecture \ref{conj:main} is well-posed.

\subsection{Twisted Brauer quotients}\label{subsec:twisted_Brauer_quotients}  Let $\tS$ be a  $p$-group and $\cA$ a bifree $\tS$-algebra with $(\tS,\tS)$-invariant $\cO$-basis $Y$.

We wish to understand what properties of the $(\tS,\tS)$-biset $Y$ can be derived from the $(\cO\tS,\cO\tS)$-bimodule structure of $\cA$, and more generally the interior $\tS$-algebra structure of $\cA$.  The Brauer quotient construction provides us with a key bridge between the worlds of algebras and  bisets, so it will be helpful to review some basic facts here.  The core content of this section is standard, see, e.g., \cite[\S\S 11,27]{ThevenazBook}.

Observe that viewing $\cA$ as a $(\cO\tS,\cO\tS)$-bimodules is equivalent to giving it the $\cO(\tS\times\tS)$-module structure determined by $(\ts_1,\ts_2)\odot a:=\ts_1\cdot a\cdot \ts_2^{-1}$ for $\ts_1,\ts_2\in \tS$ and $a\in \cA$.  Similarly, the $(\tS,\tS)$-biset structure of $Y$ is the same as thinking of $Y$ as an $(\tS\times\tS)$-set.  We shall  move back and forth freely between  these perspectives.

Recall that, for $\mathtt{U}\leq\tS\times\tS$, the \emph{Brauer quotient}  is the $\kk$-module
\[
\cA(\mathtt{U}):=\cA^\mathtt{U}\Big/\left(\mathfrak{m}.\cA^{\mathtt{U}}+\sum_{\mathtt{V}\lneq \mathtt{U}}\tra_\mathtt{V}^\mathtt{U}(\cA^\mathtt{V})\right),
\]
where $\cA^\mathtt{U}$ denotes the $\mathtt{U}$-fixed submodule of $\cA$, $\tra_\mathtt{V}^\mathtt{U}:\cA^\mathtt{V}\to\cA^\mathtt{U}$  the relative trace map, and $\mathfrak{m}$  the unique maximal ideal of $\cO$.  The \emph{Brauer homomorphism} of $\cA$ at $\mathtt{U}$ is defined to be the natural surjection $\br_\mathtt{U}:\cA^\mathtt{U}\to\cA(\mathtt{U})$.
 
The assumption that $\cA$ is a bipermutation $\tS$-algebra gives us a great deal of control over its Brauer quotients.

\begin{lemma}\label{lem:bipermutation_algebra_Brauer_quotient_general}
Let $\cA$ be a bipermutation  $\tS$-algebra with $(\tS,\tS)$-invariant $\cO$-basis $Y$.  Then for any $\mathtt{U}\leq \tS\times \tS$, the set $\{\br_\mathtt{U}(y)\ \big|\ y\in Y^\mathtt{U}\}$ is a $\kk$-basis for $\cA(\mathtt{U})$. 
\end{lemma}

The bifreeness of $\cA$ as an interior $\tS$-algebra tells us even more about its Brauer quotients:   By Lemma \ref{lem:bipermutation_algebra_Brauer_quotient_general}, $\cA(\mathtt{U})=0$ unless $\mathtt{U}$ is a twisted diagonal subgroup.  All of our algebras are bifree, so we henceforth  consider only Brauer quotients of twisted diagonal subgroups.  As with bifree bisets, we have some specialized notation for this case:  Given $\Delta(\varphi,\tP)\leq\tS\times\tS$, we write $^\varphi\!\!\cA^\tP$ for the fixed-point submodule $\cA^{\Delta(\varphi,\tP)}$, $\cA(\varphi)$ for the twisted Brauer quotient $\cA(\Delta(\varphi,\tP))$, and $\br_\varphi:{^\varphi\!\!\cA^\tP}\to\cA(\varphi)$ for the corresponding Brauer map.  

Observe that if $\varphi=\iota_\tP^\tS$ is the natural inclusion map, $\cA(\iota_\tP^\tS)$ is the standard Brauer quotient of $\cA$ at $\tP$, for which we  use the conventional notation $\cA(\tP)$.  The more general  $\cA(\varphi)$ will be called the \emph{($\varphi$-)twisted Brauer quotient} of $\cA$.  

In the context of testing for the stability of $Y$ relative to some fusion system, we observe the following reformulation of Lemma \ref{lem:bipermutation_algebra_Brauer_quotient_general}:

\begin{lemma}\label{lem:Brauer_quotient_dimensions}
Let $\cA$ be a bifree $\tS$-algebra with $(\tS,\tS)$-invariant $\cO$-basis $Y$.
 For any twisted diagonal subgroup $\Delta(\varphi,\tP)\leq \tS\times \tS$, we have
\[
\left|{^\varphi Y^\tP}\right|=\dim_\kk\cA(\varphi).
\]
\end{lemma}

The algebra structure of $\cA$ induces a multiplication on $\cA(\tP)$, but not on an arbitrary twisted Brauer quotient $\cA(\varphi)$.  Nevertheless, our focus on twisted diagonal subgroups does allow us to make use of the multiplication in $\cA$ to compare different twisted Brauer quotients.  If  $\tP,\tQ\leq \tS$ are two subgroups and $\varphi:\tP\hookrightarrow \tS$, $\psi:\tQ\hookrightarrow \tS$  two group injections such that $\varphi \tP\leq \tQ$, the multiplication of $\cA$ induces a bilinear pairing
\[
{^\psi\!\!\cA^\tQ}\times{^\varphi\!\!\cA^\tP}\to{^{\psi\varphi}\!\!\cA^\tP}.
\]
Moreover, this pairing is compatible with the restriction and transfer morphisms between fixed-point submodules for twisted diagonal subgroups, which compatibility yields:

\begin{lemma}\label{lem:twisted_brauer_quotient_connection}
For $\Delta(\varphi,\tP),\Delta(\psi,\tQ)\leq\tS\times\tS$ twisted diagonal subgroups  with $\tQ=\varphi \tP$,  the multiplication of $\cA$ induces $\cA(\psi)\otimes\cA(\varphi)\to\cA(\psi\varphi)$: For $a_\varphi\in{^\varphi\!\!\cA^\tP}$ and $a_\psi\in{^\psi\!\!\cA^\tQ}$ we have
\[
\br_\psi(a_\psi).\br_\varphi(a_\varphi)=\br_{\psi\varphi}(a_\psi.a_\varphi).
\]
In particular, $\cA(\varphi)$ is an $(\cA(\tQ),\cA(\tP))$-bimodule.
\end{lemma}

We  close this section with the observation that twisted Brauer quotients may be used to define something like a fusion system: The \emph{fixed-point fusion presystem} on $\tS$ induced by $\cA$, denoted $\fF_\tS(\cA)$,  is the partial category (i.e., not all  composites of morphisms need be defined) whose objects are the subgroups of $\tS$ and whose homsets are given by
\[
\Hom_{\fF_\tS(\cA)}(\tP,\tQ):=\left\{\varphi:\tP\to \tQ\ \big|\ \cA(\varphi)\neq 0\right\}.
\]
In terms of the $(\tS,\tS)$-invariant $\cO$-basis $Y$ of $\cA$, we may also write
\[
\Hom_{\fF_\tS(\cA)}(\tP,\tQ)=\left\{\varphi:\tP\to \tQ\ \big|\ {^\varphi Y^\tP}\neq\emptyset\right\},
\]
so that $\fF_\tS(\cA)$ may be seen as the ``fixed-point pre-fusion system'' of the $(\tS,\tS)$-biset $Y$, as defined in \cite[Definition 5.4(b)]{RagnarssonStancuIdempotents}.

Observe that the choice of the target $\tQ$ containing $\varphi\tP$ does not affect whether $\cA(\varphi)$ is nontrivial.  In other words, every morphism of $\fF_\tS(\cA)$ factors as an isomorphism of groups (although not necessarily an isomorphim in $\fF_\tS(\cA)$) followed by an inclusion of subgroups. We are particularly interested in the case where this factorization occurs in $\fF_\tS(\cA)$.  We say that $\cA$ is \emph{divisible} if $\fF_\tS(\cA)$ is a divisible category in the sense of Puig (cf. \cite{PuigBook}).  In other words, we require that $\fF_\tS(\cA)$ be a category that contains all subgroup inclusions and in which every morphism factors as a $\fF_\tS(\cA)$-isomorphism followed by an inclusion.  We note the following sufficient condition for $\cA$ to be divisible.

\begin{lemma}\label{lem:divisibility_sufficiency}
Let $\cA$ be a bifree  $\tS$-algebra.  If $\cA(\tS)\neq 0$ and there exists an $(\tS,\tS)$-invariant $\cO$-basis $Y\subseteq\cA^\times$,  then $\cA$ is divisible.
\begin{proof} 
The assumption that $\cA(\tS)\neq 0$ implies that $^{\id_\tS}Y^\tS=Y^\tS\neq\emptyset$ by Lemma \ref{lem:Brauer_quotient_dimensions}, and as $Y^\tS\subseteq Y^\tP$ for all $\tP\leq \tS$, it follows that all inclusions $\iota_\tP^\tS:\tP\leq \tS$, and indeed all inclusions between subgroups of $\tS$, are morphisms of  $\fF_\tS(\cA)$. 

As $\cA(\tP)\neq 0$ for all $\tP\leq \tS$, we have that every unit of $\cA$ has nonzero image in all twisted Brauer quotients to which it maps, i.e., if $u\in\cA^\times\cap{^\varphi\!\!\cA^\tP}$, then $0\neq\br_\varphi(u)\in\cA(\varphi)$.   To see this, note that $u^{-1}\in{^{\varphi^{-1}}\!\!\!\cA^{\varphi \tP}}$, so the multiplicative pairing of twisted Brauer quotients (Lemma \ref{lem:twisted_brauer_quotient_connection}) implies
\[
0\neq\br_\tP(1_\cA)=\br_P(u^{-1}.u)=\br_{\varphi^{-1}}(u^{-1}).\br_\varphi(u)
\]
as $\br_\tP(\cA^\tP)=\cA(\tP)\neq 0$.   

For $\varphi\in\Hom_{\fF_\tS(\cA)}(\tP,\tQ)$ such that $\varphi\tP=\tQ$, we have by definition $\cA(\varphi)\neq 0$, and therefore $\varphi$ is realized by some $y\in {^\varphi Y^\tP}$  (Lemma \ref{lem:Brauer_quotient_dimensions}).  The inverse $y^{-1}\in\cA^\times\cap{^{\varphi^{-1}}\!\!\!\cA^{\tQ}}$ has nonzero image in $\cA(\varphi^{-1})$, so $\varphi^{-1}\in\Hom_{\fF_\tS(\cA)}(\tQ,\tP)$.  It  follows that any morphism of $\fF_\tS(\cA)$ that is an isomorphism of groups  is also an isomorphism in the fixed-point fusion presystem.

Finally, we verify that all compositions in $\fF_\tS(\cA)$ are defined.  This uses the same argument as above:  If $\varphi\in\Hom_{\fF_\tS(\cA)}(\tP,\tQ)$ is realized by $y_\varphi\in\cA^\times\cap{^{\varphi}\!\!\cA^\tP}$ and $\psi\in\Hom_{\fF_\tS(\cA)}(\tQ,\tR)$ is realized by $y_\psi\in\cA^\times\cap{^\psi\!\!\cA^\tQ}$, the product $y_\psi.y_\varphi\in\cA^\times\cap{^{\psi\varphi}\!\!\cA^\tP}$, being a unit, has nonzero image in $\cA(\psi\varphi)$, and thus $\psi\varphi\in\Hom_{\fF_\tS(\cA)}(\tP,\tR)$.  This completes the proof.
\end{proof}
\end{lemma}

Returning to our true object of interest:  Reduction modulo the maximal ideal of $\cO$ yields the $\kk$-algebra $\overline\cS$, which is again an interior $\tD$-algebra.  If $X$ is a $(\tD,\tD)$-invariant $\cO$-basis of $\cS$, then its image $\overline X$ is a $(\tD,\tD)$-invariant $\kk$-basis of $\overline\cS$, which is isomorphic to $X$ as a $(\tD,\tD)$-biset.  It is  easy to check that $\overline\cS$ is a symmetric bifree $\tD$-algebra with $\overline\cS(\tD)\neq 0$.  It follows from \cite[Proposition 6]{GelvinCharacteristicBlockBasisOnline} that $\overline X$, and hence $X$, is symmetric as a $(\tD,\tD)$-biset, from which we conclude that every morphism of $\fF_\tD(\cS)$ factors as an isomorphism followed by an inclusion.  It is not immediately clear that composition of morphisms of $\fF_\tD(\cS)$ is always defined, so we are not yet able to say that $\cS$ is a divisible category.  If we were able to show the existence of a unital $(\tD,\tD)$-invariant $\cO$-basis, the divisibility of the source algebra would follow from Lemma \ref{lem:divisibility_sufficiency}, but of course this would prove Conjecture \ref{conj:main_strong_i} directly.  Luckily, we can realize the divisibility of $\cS$ by other means, as we will see in the next section.

\subsection{The current state of Conjecture \ref{conj:main}}\label{subsec:state_of_the_art}

Proving that a $(\tD,\tD)$-stable $\cO$-basis $X$ of  $\cS$ is $\cF_\tD(b)$-characteristic amounts to verifying the five conditions of Definition \ref{def:F_characteristic}:

\begin{enumerate}
\item {\bf $X$ is bifree:}  $\cS$ is a corner subalgebra of $\cG:=kG$, which has $G$ itself as a bifree $(\tD,\tD)$-stable basis.  Corollary \ref{cor:basis_well_defined}\ref{cor:basis_well_defined_iii} implies that $X$ is isomorphic to a $(\tD,\tD)$-subbiset of $\tG$, and hence is bifree.
\item {\bf $X$ is symmetric:}  This was proved in the final paragraph of \S\ref{subsec:twisted_Brauer_quotients}.
\item {\bf $X$ is $\cF_\tD(b)$-generated:}  This will be proved in Corollary \ref{cor:generation_of_source_basis}.
\item {\bf $X$ is $\cF_\tD(b)$-stable:}  This is open.  

The upcoming Theorem \ref{thm:source_algebras_are_divisible} may be seen as step in the right direction:  From the fixed-point formulation of $\cF_\tD(b)$-stability,  for all $\varphi\in\Hom_{\cF_\tD(b)}(\tP,\tQ)$ we must have $\left|{^\varphi\! X^P}\right|=\left|X^P\right|$.  As the latter term is nonzero for all $\tP\leq \tD$, the least we could hope for is that the former is as well.  This is a direct consequence of the equality $\fF_\tD(\cS)=\cF_\tD(b)$, which says that for our morphism $\varphi$ there is some $x\in X$ such that $\Stab(x)\geq\Delta(\varphi,\tP)$, and hence $^\varphi\!X^\tP\neq\emptyset$.  
\item {\bf $|X|/|\tD|$ is prime to $p$:}  This is due to Puig \cite[Proposition 14.6]{PuigConstructionOfModules}, and can also be found in \cite[Corollary 44.8]{ThevenazBook}.

In fact, more is true:  The elements of $X$ with point-stabilizers of the form $\Delta(\alpha,\tD)$ for some $\alpha\in\Aut_{\cF_\tD(b)}(\tD)$ can be explicitly computed in terms of the block fusion system, cf. \cite[Theorem 44.3(a)]{ThevenazBook}.  One observes that the $(\tD,\tD)$-orbits that appear are precisely the `top level' of the minimal characteristic biset $\Omega_{\cF_\tD(b)}$ of \cite{GelvinReehMinimalCharacteristicBisets}, and that each orbit occurs with multiplicity one.  If follows that, if Conjecture \ref{conj:main} were true, the $(\tD,\tD)$-stable $\cO$-basis of a source algebra would consist of exactly one copy of the minimal characteristic biset of the block fusion system, plus some unknown number of `lower level' terms.
\end{enumerate}

While the above state of the art is essentially unchanged since Linckelmann and Webb's exploration of characteristic bisets circa 2000, we should note that Conjecture \ref{conj:main} can  be verified in several extremal situations.  In fact, all of the following are examples where Conjecture \ref{conj:main_strong_i} can be shown to hold, if the forward referencing to Proposition \ref{prop:unit_stabilizers_implies_unit_basis} and Corollary \ref{cor:strong_1_implies_main} will be forgiven.  

\begin{enumerate}
\item If $\cS=\cB$, so that the source algebra is as large as possible, the $\cF_\tD(b)$-stability of $\cS$ is a consequence of $\cB$'s being an interior $\tG$-algebra.
\item If $\textrm{rank}_\cO\phantom{.}\cB=|\tG|_p^2/|\tD|$, so that the source algebra is as small as possible,  it is shown in \cite{LinckelmannBlocksOfMinimalDimension}  that $\cS\cong\cO\tS$, where the result is clear. 
\item  If $\tD\trianglelefteq\tG$, the defect group is as small as possible, as we always have $O_p(\tG)\leq\tD$.  In this case it is known (cf. \cite[Theorem 45.12]{ThevenazBook}) that $\cS$ is isomorphic to a twisted group algebra.  The $(\tD,\tD)$-invariant $\cO$-basis may therefore be taken to be elements of the underlying group, which are units in $\cS$.
\item  If $\cF_\tD(b)=\cF_\tD(\tD)$, the block fusion system is nilpotent and hence as small as possible. As $\cF_\tD(\tD)$-stability is an  empty condition, $\cF_\tD(b)$-stability is automatic.
\end{enumerate}

More generally, proving Conjecture \ref{conj:main} boils down to establishing the $\cF_\tD(b)$-stability of the $(\tD,\tD)$-invariant $\cO$-basis of $\cS$.  Unfortunately, we do not currently see any way to verify this property outside of the above contexts, so in the next section we will identify additional structural properties that would imply it.  

\section{Block fusion systems and local Puig categories}\label{sec:Block_fusion_and_local_categories}

Now that we have recalled the definition of a characteristic biset for a fusion system, and how to connect this notion with the structure of interior algebras, we should name the particular fusion system of interest.  We should emphasize  that, aside from proving that source algebras of blocks are divisible in Theorem \ref{thm:source_algebras_are_divisible}, the content of this section is already established in the literature and can be found in \cite{ThevenazBook} or \cite{LinckelmannBookI, LinckelmannBookII}. However, our perspective is firmly rooted in an inital choice of source idempotent, so the concept of ``exomorphism'' that informs so much of Puig's work is not particularly useful to us.  We therefore give a thorough treatment of the following established material to demonstrate that the choices we've made can be lived with, as well as to give a reference to the technical facts that will be essential in the final sections of this paper.

Given a block $b$ of $\cG=\cO\tG$ with defect group $\tD$, the block fusion system should be a fusion system on $\tD$ whose morphisms `respect $b$,' in some sense.  The right way to formulate this is by taking as the objects of the block fusion system not the subgroups of $\tD$, but rather subgroups with an additional datum that more directly relates to $b$.  The resulting structure is a $\tG$-poset via the conjugation action of $\tG$, and the morphisms of the block fusion system are then defined to be the $\tG$-conjugations that respect the partial order.  We recall these facts below.

Let $\cA$ be an  $\tG$-algebra.  A \emph{Brauer pair} of $\cA$ is a pair $(\tP,e_\tP)$, where $\tP\leq \tG$ is a subgroup and $e_\tP$  a block idempotent of $\cA(\tP)$.  As $\cA(\tP)=0$ unless $\tP$ is a $p$-group, we will assume this condition without further comment.  When $\cA=\cG$ is a group algebra, an easy consequence of Lemma \ref{lem:Brauer_quotient_dimensions} is that the Brauer morphism $\br_\tP:\cG^\tP\to\cG(\tP)$ induces an isomorphism of $N_\tG(\tP)$-algebras $\kk C_\tG(\tP)\cong\cG(\tP)$, so that a Brauer pair may be thought of as a $p$-subgroup of $\tG$ and a block of that subgroup's centralizer.  This was the view originally taken by Alperin and Brou\'e \cite{AlperinBroueBlockFusion}, who then defined a partial order on the set of Brauer pairs using the $N_\tG(\tP)$-algebra structure of $\kk C_\tG(\tP)$ to relate the blocks of $\kk C_\tG(\tP)$ and $\kk C_\tG(\tQ)$ when $\tQ\trianglelefteq \tP$.  Brou\'e and Puig \cite{BrouePuigLocalStructure} then observed that only the $p$-permutation $\tG$-algebra structure of $\cG$ is needed to define the partial order.  This second perspective makes use of notions of central importance to this paper, so it is the one we adopt.

For an arbitrary $\cO$-free $\cO$-algebra $\widetilde\cA$, a \emph{point} of $\widetilde\cA$ is an $\widetilde\cA^\times$-conjugacy class of primitive idempotents in $\widetilde\cA$.  If $i\in\widetilde\cA$ is a primitive idempotent, we  write $[i]$ for the point containing $i$.  We also write $\cP(\widetilde\cA)$ for the set of points of $\widetilde\cA$.

The $\tG$-algebra structure on $\cA$ allows us to speak of a point for each subgroup of $\tG$ and each fixed-point subalgebra of $\cA$.   A \emph{pointed group} on $\cA$ is a pair $\tH_\beta:=(\tH,\beta)$ of a subgroup $\tH\leq \tG$ and a point $\beta$ of $\cA^\tH$.  If $\tP\leq\tG$ is a $p$-subgroup, the point $\gamma\in\cA^\tP$ is \emph{local} if $\br_\tP(\gamma)\neq0$, or equivalently if $\br_\tP(i_\tP)\neq 0$ for any $i_\tP\in\gamma$.  We write $\cL\cP(\cA^\tP)$ for the set of local points of $\cA^\tP$. A \emph{local pointed group} on $\cA$ is a pointed group $\tP_\gamma$ on  $\cA$ whose point is local.

The set of local pointed groups on $\cA$ is finite, and may be thought of as a generalization of the set of $p$-subgroups of $\tG$.  It also has a natural partial order:  $\tQ_\delta\leq \tP_\gamma$ means $\tQ\leq \tP$ and there exist $j_\tQ\in\delta$ and $i_\tP\in\gamma$ such that $j_\tQ\leq i_\tP$, i.e., $j_\tQ=j_\tQ.i_\tP.j_\tQ$.  By analogy with the Brown poset of $p$-subgroups of $\tG$, we write $s(\tG,\cA)$ for the poset of local pointed groups on $\cA$, and  refer to this as the \emph{pointed Brown poset}.  

The poset structure on $s(\tG,\cA)$ induces the partial order on the set of Brauer pairs of $\cA$.  As $\gamma\in\cL\cP(\cA^\tP)$ consists of primitive idempotents, and $\br_\tP(\gamma)\neq 0$, there is a unique block $e_\gamma$ of $\cA(\tP)$ such that $\br_\tP(i_\tP)\leq e_\gamma$ for some (hence all) $i_\tP\in\gamma$.  We then say that the local pointed group $\tP_\gamma$ is \emph{associated to} the Brauer pair $(\tP,e_\gamma)$. For Brauer pairs $(\tQ,e_\tQ)$ and $(\tP,e_\tP)$ on $\cA$, we write $(\tQ,e_\tQ)\leq (\tP,e_\tP)$ if there exist $\tQ_\delta,\tP_\gamma\in s(\tG,\cA)$ with $\tQ_\delta$ associated to $(\tQ,e_\tQ)$, $\tP_\gamma$ associated to $(\tP,e_\tP)$, and $\tQ_\delta\leq\tP_\gamma$.

The group $\tG$ acts on the pointed Brown poset $s(\tG,\cA)$ by conjugation.   It follows that the relation $\leq$ on the set of Brauer pairs of $\cA$ is $\tG$-equivariant:  If $(\tQ,e_\tQ)\leq (\tP,e_\tP)$ and $\tg\in \tG$, then $^\tg(\tQ,e_\tQ)\leq{^\tg(\tP,e_\tP)}$.  The key properties of the relation $\leq$  are (cf. \cite[\S40]{ThevenazBook}):

\begin{prop}\label{prop:Brauer_partial_order_properties}
Let $\cA$ be a $p$-permutation $\tG$-algebra.
\begin{enumerate}
\item\label{prop:Brauer_partial_order_properties_i}  The relation $\leq$ is a $\tG$-equivariant partial order on the set of Brauer pairs of $\cA$.
\item\label{prop:Brauer_partial_order_properties_ii}  $(\tQ,e_\tQ)\leq(\tP,e_\tP)$ if and only $\tQ\leq \tP$ and there is some primitive idempotent $i_\tP\in\cA^\tP$ such that $\br_\tP(i_\tP)\leq e_\tP$ and $\br_\tQ(i_\tP)\leq e_\tQ$.
\item\label{prop:Brauer_partial_order_properties_iii}  Given a Brauer pair $(\tP,e_\tP)$ of $\cA$ and $\tQ\leq \tP$, there is a \emph{unique} block $e_\tQ$ of $\cA(\tQ)$ such that $(\tQ,e_\tQ)\leq(\tP,e_\tP)$.
\item\label{prop:Brauer_partial_order_properties_iv}  All maximal elements in a connected component in the poset of Brauer pairs of $\cA$ are $\tG$-conjugate.
\end{enumerate}
\end{prop}

Returning to the case of the group algebra $\cG$, as the Brauer quotient $\cG(1)=\cG$ we have that the blocks of $\cG$ are precisely the second coordinates of Brauer pairs whose subgroups are trivial.  If $(\tP,e_\tP)$ is an arbitrary Brauer pair of $\cG$, Proposition \ref{prop:Brauer_partial_order_properties}\ref{prop:Brauer_partial_order_properties_iii} implies that there is a unique block $b$ of $\cG$ such that $(1, b)\leq (\tP,e_\tP)$.  It follows that the connected components of the poset of Brauer pairs of $\cG$ are parameterized by the blocks of $\cG$.  Moreover, as $^\tg(1, b)=(1, b)$ for all $\tg\in \tG$, each connected component is closed under $\tG$-conjugation, and hence Proposition \ref{prop:Brauer_partial_order_properties}\ref{prop:Brauer_partial_order_properties_iv} yields that the set of maximal elements of a component  is the entire $\tG$-orbit of some Brauer pair $(\tD,e_\tD)$.  

A \emph{defect group} of the block $b$ is  any $\tD\leq \tG$ such that there is a maximal Brauer pair $(\tD,e_\tD)$ with $(1, b)\leq(\tD,e_\tD)$.  Once such a maximal Brauer pair is fixed, Proposition \ref{prop:Brauer_partial_order_properties}\ref{prop:Brauer_partial_order_properties_iii} gives for each $\tP\leq \tD$ a unique block $e_\tP$ of $\cG(\tP)$ such that $$(1, b)\leq(\tP,e_\tP)\leq(\tD,e_\tD).$$  In particular, we may identify the interval $[(1, b),(\tD,e_\tD)]$ in the Brauer pair poset with the subgroup poset of $\tD$ itself.  The \emph{block fusion system} of $b$ is then defined to be the fusion system $\cF_\tD(b)$ on $\tD$ whose morphisms are given by
\[
\Hom_{\cF_\tD(b)}(\tP,\tQ)=\left\{\varphi:\tP\to \tQ\ \big|\ \exists\ \tg\in \tG:\ {^\tg(\tP,e_\tP)\leq (\tQ,e_\tQ)},\ \varphi=c_\tg|_\tP\right\}.
\]
As all maximal Brauer pairs lying above $(1, b)$ are $G$-conjugate, it is easy to see that $\cF_\tD(b)$ is well-defined up to isomorphism of fusion systems, independent of the choice of $(\tD,e_\tD)$.  (However, we remind the reader that a \emph{particular} choice of block fusion system does depend on the defect Brauer pair, although we suppress this dependence in our notation.)  Moreover, $\cF_\tD(b)$ is saturated (cf. \cite[Part IV, Theorem 3.2]{AKO}),  so we may make use of Alperin's fusion theorem \cite{AlperinFusion} and other such structural results.  In particular, there exist $\cF_\tD(b)$-characteristic bisets, and there is hope for Conjecture \ref{conj:main}.

We use the local pointed group definition of the partial order on the set of Brauer pairs because most of our results will rely on analyses of various primitive idempotent decompositions of fixed-point subalgebras.  The local pointed groups can be organized into  \emph{local categories} whose structures provide the appropriate framework for comparing idempotent decompositions.  We describe these categories now.

Let $\tG$ be a finite group, $\tS\in\Syl_p(\tG)$ a Sylow $p$-subgroup, and $\cA$ a  $\tp$-permutation interior  $\tG$-algebra.  In particular, $\cA$ is  an  interior permutation $\tS$-algebra. The \emph{local  $\tG$-Puig category} $\fL_\tS^\tG(\cA)$ has as its objects the pointed Brown poset $s(\tS,\cA)$ and homsets defined by
\[
\Hom_{\fL_\tS^\tG(\cA)}(\tP_\gamma,\tQ_\delta):=\left\{\varphi:\tP\to \tQ\ \big|\ \exists\ \tg\in \tG:\ {^\tg(\tP_\gamma)\leq \tQ_\delta},\ \varphi=c_\tg|_\tP\right\}.
\]
Morphisms of $\fL_\tS^\tG(\cA)$ are called \emph{$\tG$-fusions in $\cA$}.

To emphasize the connection with  fusion systems, note that if $\cA=\kk$ is the trivial $\tG$-algebra, then every $p$-subgroup of $\tG$ has a unique local point on $\cA$ and one can identify $\fL_\tS^\tG(\kk)$ with  $\cF_\tS(\tG)$.

While both  group  and block algebras are interior $\tG$-algebras, the best we can say of the source algebra is that it is an interior algebra for its defect group.  We need a variant of the local  $\tG$-Puig category that recovers the part of $\fL_\tS^\tG(\cG)$ associated with $b$, namely $\fL_\tD^\tG(\cB)$, but that is still sensible when applied to $\cS$.  Such a construction is provided below.

Let $\tS$ be a $p$-group and $\cA$ an interior $\tS$-algebra.  The \emph{local unital Puig category} $\fL_\tS^\times(\cA)$ is the category with objects $s(\tS,\cA)$  and whose morphisms are defined as follows:  Given local pointed groups $\tP_\gamma,\tQ_\delta\in s(\tS,\cA)$, let $\Hom_{\fL_\tS^\times(\cA)}(\tP_\gamma,\tQ_\delta)$ be the set of injective group maps $\varphi:\tP\hookrightarrow \tQ$ such that for some $i_\tP\in\gamma$ and $j_\tQ\in\delta$ there is a unit $u\in\cA^\times$ satisfying
\begin{enumerate}
\item $^ui_\tP\in\cA^{\varphi \tP}$,
\item $^ui_\tP\leq j_\tQ$, and
\item $^u({\tp\cdot i_\tP})=\varphi(\tp)\cdot {^u i_\tP}$ for all $\tp\in \tP$.
\end{enumerate}
Morphisms of $\fL_\tS^\times(\cA)$ are called \emph{$\cA^\times$-fusions}. One readily verifies that  whether  $\varphi:\tP\hookrightarrow \tQ$ is an $\cA^\times$-fusion is independent of the choice of idempotents $i_\tP\in\gamma$ and $j_\tQ\in\delta$.  

It is   easy to see that all inclusions $\tP_\gamma\leq \tQ_\delta$ of local pointed groups are $\cA^\times$-fusions.  Similarly, if $\varphi:\tP\hookrightarrow\tQ$ is both an $\cA^\times$-fusion and a group isomorphism,  then $\varphi$ is an isomorphism in $\fL_\tS^\times(\cA)$.   The same comments apply to $\fL_\tS^\tG(\cA)$, and  so local Puig categories are divisible. Isomorphisms in local Puig categories will be referred to as \emph{isofusions}.

The divisibility of local Puig categories allows us to concentrate on isofusions to simplify the proofs of many results, not the least because the $\cA^\times$-isofusions  have  an easier characterization. A group isomorphism $\varphi:\tP\xrightarrow\cong \tQ$ between subgroups of $\tS$ lies in $\Iso_{\fL_\tS^\times(\cA)}(\tP_\gamma,\tQ_\delta)$ if and only if for some (hence all) $i_\tP\in\gamma$ and $j_\tQ\in\delta$, there exists  $u\in\cA^\times$ such that 
 \[
 ^u(\tp\cdot i_\tP)=\varphi(\tp)\cdot j_\tQ\textrm{ for all }\tp\in \tP.
 \]
In particular, taking $\tp=1$, we have $^ui_\tP=j_\tQ$.  Thus $\cA^\times$-isofusions can be viewed as a refinement of the conjugacy relation on idempotents.  As idempotents in a finite dimensional $\cO$-algebra are conjugate if and only if they are associate, it is unsurprising that we have the following  reformulation:

\begin{lemma}\label{lem:Puig_morphisms_via_transpotents}
Let $\cA$ be an interior $\tS$-algebra and $\tP_\gamma,\tQ_\delta\in s(\tS,\cA)$ local pointed subgroups with representative idempotents $i_\tP\in\gamma$ and $j_\tQ\in\delta$. Any group isomorphism $\varphi:\tP\xrightarrow{\cong} \tQ$ lies in $\Iso_{\fL_\tS^\times(\cA)}(\tP_\gamma,\tQ_\delta)$ if and only if  there exist
\[
s_\varphi\in j_\tQ.{^\varphi\!\!\cA^\tP}.i_\tP\qquad\textrm{and}\qquad
t_{\varphi^{-1}}\in i_P.{^{\varphi^{-1}}\!\!\!\cA^\tQ}.j_\tQ
\]
such that $i_\tP=t_{\varphi^{-1}}.s_\varphi$ and $j_\tQ=s_\varphi.t_{\varphi^{-1}}$.
\begin{proof}
First suppose that $\varphi:\tP_\gamma\xrightarrow\cong\tQ_\delta$ is an $\cA^\times$-isofusion and let $u\in\cA^\times$ realize this fact.  We have $^ui_\tP=j_\tQ$, or $u.i_\tP=j_\tQ.u$, so set
\[
s_\varphi:=j_\tQ.u.i_\tP=j_\tQ.u=u.i_\tP\qquad\textrm{and}\qquad
t_{\varphi^{-1}}:=i_\tP.u^{-1}.j_\tQ=i_\tP.u^{-1}=u^{-1}.j_\tQ
\]
to obtain elements lying in $j_\tQ.\cA.i_\tP$ and $i_\tP.\cA.j_\tQ$, respectively.  These elements also live in the desired fixed-point submodules of $\cA$:  We have $^u(\tp\cdot i_\tP)=\varphi(\tp)\cdot j_\tQ$ for all $\tp\in \tP$, so  
\[
s_\varphi\cdot \tp=(u.i_\tP)\cdot \tp=u.(\tp\cdot i_\tP)={^u(\tp\cdot i_\tP)}\cdot u=(\varphi(\tp)\cdot j_\tQ).u=\varphi(\tp)\cdot s_\varphi,
\]
and a similar computation yields $t_{\varphi^{-1}}\cdot \tq=\varphi^{-1}(\tq)\cdot t_{\varphi^{-1}}$ for all $\tq\in \tQ$. We then have
\[
t_{\varphi^{-1}}.s_\varphi=(i_\tP.u^{-1}).(u.i_\tP)=i_\tP\qquad\textrm{and}\qquad
s_\varphi.t_{\varphi^{-1}}=(j_\tQ.u).(u^{-1}.j_\tQ)=j_\tQ,
\]
proving the ``only if'' implication.

Conversely, suppose we are given elements $s_\varphi$ and $t_{\varphi^{-1}}$ as prescribed.  As $i_\tP=t_{\varphi^{-1}}.s_\varphi$ and $j_\tQ=s_\varphi.t_{\varphi^{-1}}$, the idempotents $i_\tP$ and $j_\tQ$ are associate in $\cA$.  An equivalent formultation of $i_\tP$'s and $j_\tQ$'s being associate is that there is an isomorphism $\cA.i_\tP\cong \cA.j_\tQ$ of left $\cA$-modules:  If $f:\cA.i_\tP\to\cA.j_\tQ$ is  a left $\cA$-module isomorphism it is easy to check that $i_\tP=f(i_\tP).f^{-1}(j_\tQ)$ and $j_\tQ=f^{-1}(j_\tQ).f(i_\tP)$, while the factorizations $i_\tP=t.s$ and $j_\tQ=s.t$ give rise to inverse isomorphisms $\cA.i_\tP\to \cA.j_\tQ:a\mapsto a. t$ and $\cA.j_\tQ\to \cA.i_\tP:a\mapsto a.s$.  

Let $\tilde \imath_\tP:=1_\cA-i_\tP$ and $\tilde \jmath_\tQ:=1_\cA-j_\tQ$ be the complementary idempotents to $i_\tP$ and $j_\tQ$.  As $\cA.i_\tP\cong \cA.j_\tQ$ as left $\cA$-modules and $\cA.i_\tP\oplus\cA.\tilde\imath_\tP=\cA=\cA.j_\tQ\oplus\cA.\tilde\jmath_\tQ$, the Krull-Schmidt theorem implies that $\cA.\tilde \imath_\tP\cong\cA.\tilde \jmath_\tQ$ as left $\cA$-modules as well.  Thus there exist $\tilde s,\tilde t\in\cA$ such that $\tilde \imath_\tP=\tilde t.\tilde s$ and $\tilde \jmath_\tQ=\tilde s.\tilde t$, although note that we say nothing about these lying in any fixed-point submodule of $\cA$.  We may, however, assume that $\tilde s\in \tilde \jmath_\tQ.\cA.\tilde \imath_\tP$ and $\tilde  t\in\tilde \imath_\tP.\cA.\tilde \jmath_\tQ$, as replacing arbitrary $\tilde s$ and $\tilde t$ with $\tilde \jmath_\tQ.\tilde s.\tilde \imath_\tP$ and $\tilde \imath_\tP.\tilde t.\tilde \jmath_\tQ$ will  yield elements whose products are still $\tilde \imath_\tP$ and $\tilde \jmath_\tQ$, depending on the order.

Set $u:=s_\varphi+\tilde s$ and $v:=t_{\varphi^{-1}}+\tilde t$.  We compute
\[
u.v=(s_\varphi+\tilde s).(t_{\varphi^{-1}}+\tilde t)=s_\varphi.t_{\varphi^{-1}}+s_\varphi.\tilde t+\tilde s.t_{\varphi^{-1}}+\tilde s.\tilde t=
j_\tQ+0+0+\tilde \jmath_\tQ=1_\cA,
\]
where the middle terms vanish as $s_\varphi.\tilde t=(s_\varphi.i_\tP).(\tilde \imath_\tP.\tilde  t)$, $\tilde s.t_{\varphi^{-1}}=(\tilde s.\tilde \imath_\tP).(i_\tP.t_{\varphi^{-1}})$, and the idempotents $i_\tP$ and $\tilde \imath_\tP$ are by construction orthogonal.  We similarly have $v.u=1_\cA$, so that $u\in\cA^\times$ and $v=u^{-1}$.  Finally, for $\tp\in \tP$, we have
\[
\begin{split}
^u(\tp\cdot i_\tP)&=(s_\varphi+\tilde  s).(\tp\cdot i_\tP).(t_{\varphi^{-1}}+\tilde  t)\\
&=s_\varphi\cdot \tp\cdot i_\tP.t_{\varphi^{-1}}+s_\varphi\cdot \tp\cdot i_\tP.\tilde t+\tilde s\cdot \tp\cdot i_\tP.t_{\varphi^{-1}}+\tilde s\cdot \tp\cdot i_\tP.\tilde t\\
&=\varphi(\tp)\cdot s_\varphi.i_\tP.t_{\varphi^{-1}}
+\varphi(\tp)\cdot s_\varphi.i_\tP.\tilde t+(\tilde s.i_\tP)\cdot \tp\cdot t_{\varphi^{-1}}
+(\tilde s.i_\tP)\cdot \tp\cdot \tilde t\\
&=\varphi(\tp)\cdot s_\varphi.t_{\varphi^{-1}}.s_\varphi.t_{\varphi^{-1}}+0+0+0\\
&=\varphi(\tp)\cdot j_\tQ^2=\varphi(\tP)\cdot j_\tQ.
\end{split}
\]
The last three terms in the third line vanish because $i_\tP.\tilde t=i_\tP.\tilde \imath_\tP.\tilde t=0$ and $\tilde s.i_\tP=\tilde s.\tilde \imath_\tP.i_\tP=0$ by our choice that $\tilde s\in\tilde \jmath_\tQ.\cA.\tilde \imath_\tP$ and $\tilde t\in\tilde\imath_\tP.\cA.\tilde \jmath_\tQ$.  Thus the unit $u$ realizes $\varphi:\tP_\gamma\xrightarrow\cong\tQ_\delta$ as an $\cA^\times$-isofusion in $\Iso_{\fL_\tS^\times(\cA)}(\tP_\gamma,\tQ_\delta)$, and we are done.
\end{proof}
\end{lemma}

It is worth noting that $\cA^\times$-isofusions is are determined by their source and the underlying group isomorphism.  (Dually, the divisibility of $\fL_\tS^\times(\cA)$ implies that the target and group map serve just as well, but we will have minimal need for this fact.)  

\begin{lemma}\label{lem:Puig_category_unique_target}
Let  $\cA$ be an interior $\tS$-algebra, $\tP,\tQ\leq \tS$ two subgroups, and $\varphi:\tP\xrightarrow\cong \tQ$ a group isomorphism.  If $\gamma\in\cL\cP(\cA^\tP)$ and $\delta,\delta'\in\cL\cP(\cA^\tQ)$ are local points such that we have both $\varphi\in\Iso_{\fL_S^\times(\cA)}(\tP_\gamma,\tQ_\delta)$ and $\varphi\in\Iso_{\fL_\tS^\times(\cA)}(\tP_\gamma,\tQ_{\delta'})$, then $\delta=\delta'$.
\begin{proof}
By Lemma \ref{lem:Puig_morphisms_via_transpotents}, we may choose $i_\tP\in\gamma$, $j_\tQ\in\delta$, and $j_\tQ'\in\delta'$, together with
\begin{eqnarray*}
s_\varphi\in j_\tQ.{^\varphi\!\!\cA^\tP}.i_\tP,&&t_{\varphi^{-1}}\in i_\tP.{^{\varphi^{-1}}\!\!\!\cA^\tQ}.j_\tQ,\\
s_\varphi'\in j_\tQ'.{^\varphi\!\!\cA^\tP} .i_P,&\textrm{and }&
t_{\varphi^{-1}}'\in i_\tP.{^{\varphi^{-1}}\!\!\!\cA^\tQ}.j_\tQ'
\end{eqnarray*}
such that $i_\tP=t_{\varphi^{-1}}.s_\varphi=t_{\varphi^{-1}}'.s_\varphi'$, $j_\tQ=s_\varphi.t_{\varphi^{-1}}$, and $j_\tQ'=s_\varphi'.t_{\varphi^{-1}}'$.  We then compute
\[
\begin{split}
(s_\varphi.t_{\varphi^{-1}}').(s_\varphi'.t_{\varphi^{-1}})=s_\varphi.i_\tP.t_{\varphi^{-1}}=s_\varphi.t_{\varphi^{-1}}=j_\tQ,\\
(s_\varphi'.t_{\varphi^{-1}}).(s_\varphi.t_{\varphi^{-1}}')=s_\varphi'.i_\tP.t_{\varphi^{-1}}'=s_\varphi'.t_{\varphi^{-1}}'=j_\tQ'.
\end{split}
\]
It is immediate that both $s_\varphi.t_{\varphi^{-1}}'$ and $s_\varphi'.t_{\varphi^{-1}}$ lie in $\cA^\tQ$.  Therefore $j_\tQ$ and $j_\tQ'$ are associate, and hence conjugate, elements of $\cA^\tQ$, and we have proved $\delta=\delta'$.
\end{proof}
\end{lemma}

Another immediate consequence of Lemma \ref{lem:Puig_morphisms_via_transpotents} is that local unital Puig categories behave well with respect to taking corner subalgebras.

\begin{prop}\label{prop:Puig_categories_of_corner_algebras_are_full}
Let  $\widehat{\cA}$ be an interior $\tS$-algebra, $\ell\in\cA^\tS$ an  idempotent, and $\cA:=\ell.\widehat{\cA}.\ell$ the interior $\tS$-algebra obtained by cutting $\widehat{\cA}$ by $\ell$.  Then the embedding $\cA\subseteq\widehat{\cA}$ induces a fully faithful functor $\fL_\tS^\times(\cA)\to\fL_\tS^\times(\widehat{\cA})$.
\begin{proof}
First observe that, as $\ell$ is fixed by $\tS$, for any $\tP\leq \tS$ we have $\cA^\tP=\ell.\widehat{\cA}^\tP.\ell$.  It follows that any idempotent lying in $\cA^\tP$ is primitive there if and only if it is primitive in $\widehat{\cA}^\tP$, and the same holds for the local condition.

Two idempotents of $\cA$ are associate in $\cA$ if and only if they are associate in $\widehat\cA$, hence any local point $\gamma\in\cL\cP(\cA^\tP)$ is contained in a single point $\widehat\gamma\in\cL\cP(\cA^\tP)$, namely the $\widehat\cA^\times$-conjugacy orbit of any $i_\tP\in\gamma$.  Similarly, if $\widehat\gamma\in\cL\cP(\widehat\cA^\tP)$ has nonempty intersection with $\cA$, then $\gamma:=\widehat\gamma\cap\cA$ is a local point in $\cL\cP(\cA^\tP)$.  Thus we have a natual bijection 
\[
\cL\cP(\cA^\tP)\leftrightarrow\left\{\widehat{\gamma}\in\cL\cP(\widehat{\cA}^\tP)\ \big|\ \widehat{\gamma}\cap\cA\neq\emptyset\right\}:\gamma\leftrightarrow\widehat{\gamma}.
\]

The functor $\fL_\tS^\times(\cA)\to\fL_\tS^\times(\widehat{\cA})$ is  defined on objects by the assignment $P_\gamma\mapsto P_{\widehat{\gamma}}$. That this assignment is fully faithful is then to say that $\varphi:\tP_\gamma\to\tQ_\delta$ is an $\cA^\times$-fusion if and only if it is an $\widehat{\cA}^\times$-fusion.  As local Puig categories are divisible, it sufffices to consider only isofusions, in which case the result is a direct application of Lemma \ref{lem:Puig_morphisms_via_transpotents}. 
\end{proof}
\end{prop}

There is another sense in which the unital local Puig category is the right object to consider for our purposes:  For any divisible $\tS$-algebra $\cA$, there is a direct connection between $\cL_\tS^\times(\cA)$ and the the fixed-point fusion presystem $\fF_\tS(\cA)$.
 
\begin{lemma}\label{lem:Puig_fusion_containment}
If $\cA$ is a divisible $\tS$-algebra with fixed-point fusion presystem $\fF_\tS(\cA)$, then for any local pointed groups $\tP_\gamma,\tQ_\delta\in s(\tS,\cA)$ we have
\[
\Hom_{\fL_\tS^\times(\cA)}(\tP_\gamma,\tQ_\delta)\subseteq\Hom_{\fF_\tS(\cA)}(\tP,\tQ).
\]
\begin{proof}
As both $\fL_\tS^\times(\cA)$ and $\fF_\tS(\cA)$ are divisible categories, it suffices to prove the result on the level of isomorphisms. Suppose that $\varphi\in\Iso_{\fL_\tS^\times(\cA)}(\tP_\gamma,\tQ_\delta)$ is witnessed by $i_\tP\in\gamma$, $j_\tQ\in\delta$, and $u\in\cA^\times$.  For $\tp\in \tP$, consider
\[
(u. i_\tP)\cdot \tp=u\cdot \tp\cdot i_\tP={^u(\tp\cdot i_\tP)}.u=\varphi(\tp)\cdot (j_\tQ.u).
\]
Since $j_\tQ={^ui_\tP}$, it follows that $u.i_\tP=j_\tQ.u$, and therefore $u.i_\tP\in{^\varphi\!\!\cA^\tP}$.  A similar computation yields $i_\tP.u^{-1}\in{^{\varphi^{-1}}\!\!\!\cA^\tQ}$.  We may therefore consider the images of these elements in their respective twisted Brauer quotients:
\[
0\neq\br_\tP(i_\tP)=\br_\tP((i_\tP.u^{-1}).(u.i_\tP))=\br_{\varphi^{-1}}(i_\tP.u^{-1}). \br_\varphi(u.i_\tP).
\]
In particular, $0\neq\br_\varphi(u.i_\tP)\in\cA(\varphi)$, so that $\varphi\in\Hom_{\fF_\tS(\cA)}(\tP,\tQ)$, as desired.
\end{proof}
\end{lemma} 

In the special case of divisible $\tS$-algebras that possess a unital $(\tS,\tS)$-invariant $\cO$-basis, we have more even precise control over  $\cA^\times$-fusions.

\begin{prop}\label{prop:local_Puig_morphisms_realized_by_unital basis}
If $\cA$ is a bifree $\tS$-algebra with unital $(\tS,\tS)$-invariant $\cO$-basis $Y$, then every $\cA^\times$-fusion is realized by an element of $Y$.
\begin{proof}
Fix $\varphi\in\Iso_{\fL_\tS^\times(\cA)}(\tP_\gamma,\tQ_\delta)$.  The existence of a unital $(\tS,\tS)$-invariant $\cO$-basis implies that $\cA$ is divisible by Lemma \ref{lem:divisibility_sufficiency}, so  $\varphi\in\Iso_{\fF_\tS(\cA)}(\tP,\tQ)$ by Lemma \ref{lem:Puig_fusion_containment}.  Thus $\cA(\varphi)\neq0$, so by Lemma \ref{lem:Brauer_quotient_dimensions} there is some $y\in{^\varphi Y^\tP}$.  Then $y^{-1}\in{^{\varphi^{-1}}\!\!\!\cA^\tQ}$, i.e., $\varphi^{-1}(\tq)\cdot y^{-1}=y^{-1}\cdot \tq$ for all $\tq\in \tQ$.  Fix $i_\tP\in\gamma$ and set $j_\tQ:={^yi_\tP}$.  Then for all $\tq\in \tQ$, we have
\[
j_\tQ\cdot \tq=y. i_\tP. y^{-1}\cdot \tq=y. i_\tP\cdot \varphi^{-1}(\tq)\cdot y^{-1}=y\cdot \varphi^{-1}(\tq)\cdot i_\tP. y^{-1}=\tq\cdot y. i_\tP. y^{-1} =\tq\cdot j_\tQ,
\]
so that $j_\tQ\in\cA^\tQ$. The local primitivity of $i_\tP$ implies the same for $j_\tQ$. From our definition of $j_\tQ$ it also follows that, for all $\tp\in \tP$, we have
\[
^y(\tp\cdot i_\tP)={^y\tp}\cdot j_\tQ=\varphi(\tp)\cdot j_\tQ,
\]
so $y$ realizes $\varphi\in\Iso_{\fL_\tS^\times(\cA)}(\tP_\gamma,\tQ_{[j_\tQ]})$.  Finally, as $\varphi\in\Iso_{\fL_\tS^\times(\cA)}(\tP_\gamma,\tQ_\delta)\cap\Iso_{\fL_\tS^\times(\cA)}(\tP_\gamma,\tQ_{[j_\tQ]})$, Lemma \ref{lem:Puig_category_unique_target} yields $\delta={[j_\tQ]}$.  Thus $y$ realizes  $\varphi$, and the claim is proved.
\end{proof}
\end{prop}

The last several results combine to prove a theorem of Puig \cite[Theorem 3.1]{PuigLocalFusions}, essentially stating that the local unital Puig category of the source algebra $\cS$ of a block contains the same data as the local $G$-Puig category of the block algebra $\cB$.

\begin{prop}\label{prop:local_Puig_categories_equivalent}
The local Puig categories $\fL_\tD^\times(\cS)$ and $\fL_\tD^G(\cB)$ are equivalent.
\begin{proof}
We first claim that the two local Puig categories of the block algebra $\cB$ are equal:  $\fL_\tD^\tG(\cB)=\fL_\tD^\times(\cB)$.  Clearly the objects are the same, so we must just compare homsets.  As both  categories are divisible, it suffices to show that $\varphi\in\Iso_{\fL_\tD^\tG(\cB)}(\tP_\gamma,\tQ_\delta)$ if and only if $\varphi\in\Iso_{\fL_\tD^\times(\cB)}(\tP_\gamma,\tQ_\delta)$.

If $\varphi:\tP_\gamma\xrightarrow\cong\tQ_\delta$ is a $\tG$-isofusion in $\cB$ realized by $\tg\in\tG$ we have $^\tg(\tP_\gamma)= \tQ_\delta$ and $\varphi=c_\tg|_\tP$.  In particular, $^\tg i_\tP\in\delta$ for all $i_\tP\in\gamma$.  We  consider the unit $u_\tg:=\tg\cdot b\in\cB$ with inverse $u_\tg^{-1}=\tg^{-1}\cdot b$, and we claim that $u_\tg$ realizes $\varphi$ as a $\cB^\times$-isofusion.  For $\tp\in \tP$ we compute
\[
^{u_\tg}(\tp\cdot i_\tP)=\tg\cdot b\cdot \tp\cdot i_\tP\cdot \tg^{-1}\cdot b={^\tg\tp}\cdot{^\tg i_\tP}\cdot b=\varphi(\tp)\cdot j_\tQ
\]
for $j_\tQ:={^\tg i_\tP}\in\delta$, which give one containment.

Suppose now that $\varphi:\tP_\gamma\xrightarrow\cong\tQ_\delta$ is a $\cB^\times$-isofusion. By Proposition \ref{prop:Puig_categories_of_corner_algebras_are_full}, $\varphi:\tP_{\widehat\gamma}\xrightarrow\cong\tQ_{\widehat{\delta}}$ is also a $\cG^\times$-isofusion, where $\tP_{\widehat{\gamma}}$ and $\tQ_{\widehat{\delta}}$ are the corresponding local pointed groups on $\cG$.   If we consider the standard $(\tD,\tD)$-invariant $\cO$-basis $\tG$ of $\cG$, the hypotheses of Proposition \ref{prop:local_Puig_morphisms_realized_by_unital basis} are satisfied, so there is some $\tg\in \tG$ realizing $\varphi$ as a $\tG$-fusion in $\cG$.  Note that the proof of that proposition further gives $\tg\in{^\varphi \tG^\tP}$, so that $\varphi=c_\tg|_\tP$.  As $b$ is central in $\cG$, we have $\widehat{\gamma}=\gamma$ and $\widehat{\delta}=\delta$.  Thus $^\tg(\tP_\gamma)={^\tg(\tP_{\widehat{\gamma}})}=\tQ_{\widehat{\delta}}=\tQ_\delta$, which completes the proof of the  claim.

It therefore suffices to prove that the inclusion $\iota:\fL_\tD^\times(\cS)\subseteq\fL_\tD^\times(\cB)$ is an equivalence of categories.  Proposition \ref{prop:Puig_categories_of_corner_algebras_are_full} already says that $\iota$ is a fully faithful embedding, so we must show it to be essentially surjective, i.e.,  every object $\tP_{\widehat\gamma}$ of $\fL_\tD^\times(\cB)$ is isomorphic to an object of the form $\tQ_{\widehat{\delta}}$ for $\tQ_\delta\in s(\tD,\cS)$ and $\widehat{\delta}\in\cL\cP(\cB^\tQ)$ the unique point containing $\delta$.  By the first claim of this proof, the term ``isomorphic'' in the previous sentence may be equally well understood to be taken in $\fL_\tD^\tG(\cB)$ as $\fL_\tD^\times(\cB)$, as the two categories are equal.  

The pointed group $\tD_{[\ell]}$ is maximal in $s(\tD,\cB)$, so by \cite[Theorem 18.3]{ThevenazBook}  $\tP_{\widehat\gamma}$ is $\tG$-subconjugate to it:  There is some $\tg\in \tG$ such that $\tQ_{\widehat\delta}:={^\tg(\tP_{\widehat\gamma})}\leq \tD_{[\ell]}$.  This containment implies that $\delta:=\widehat{\delta}\cap \cS\neq\emptyset$ is a local point of $\cS^\tQ$.  As $\tg$ clearly realizes $c_\tg:\tP_{\widehat\gamma}\xrightarrow\cong\tQ_{\widehat\delta}$ as a $\tG$-isofusion in $\cB$, this proves the essential surjectivity of $\iota$, and with it the result.
\end{proof}
\end{prop}

We are now ready to prove the main result of this section.

\begin{theorem}\label{thm:source_algebras_are_divisible}
The source algebra $\cS$ is a divisible $\tD$-algebra, and  $\fF_\tD(\cS)=\cF_\tD(b)$.
\begin{proof}
As $\cF_\tD(b)$ is a fusion system, the divisibility of $\cS$ will follow from the claimed equality. Let $X$ be a $(\tD,\tD)$-invariant $\cO$-basis of $\cS$, and take $\tG$ as a $(\tD,\tD)$-invariant $\cO$-basis of $\cG=\cO\tG$.

 $\bm{\fF_\tD(\cS)\subseteq\cF_\tD(b)}${\bf:}  This result is known in the literature, and indeed the following is essentially the argument of \cite[7.7]{LinckelmannOnSpendidDerived}, which we include for completeness. 
 
Let $x\in X$ have point-stabilizer $\Stab(x)=\Delta(\varphi,\tP)$; we will show $\varphi\in\Hom_{\cF_\tD(b)}(\tP,\tD)$.  As $\cS$ is a corner algebra of $\cG$,  $X$ is isomorphic to a $(\tD,\tD)$-subbiset of $_\tD\tG_\tD$ by Corollary \ref{cor:basis_well_defined}\ref{cor:basis_well_defined_iii}.  In particular, there is some $\tg\in \tG$ such that
$
\Stab_{(\tD,\tD)}(\tg)=\Delta(c_\tg,\tD\cap \tD^\tg)=\Stab(x).
$
The conjugation action of any such $\tg$ on $\tP$ is equal to the group map $\varphi$, so it suffices to show that $\tg$ has the desired action on Brauer pairs. Set $\tQ:={^\tg\tP}=\varphi(\tP)$, and interpret all Brauer homomorphisms as being computed in $\cG$.

Let $(\tD,e_\tD)$ be the maximal Brauer pair of $\cG$ such that $\br_\tD(\ell)\leq e_\tD$.  For  $\tP\leq \tD$, let $(\tP,e_\tP)$ be the unique Brauer pair such that $\br_\tP(\ell)\leq e_\tP$, i.e.,  $(\tP,e_\tP)\leq (\tD,e_\tD)$.  We must show that
$
^\tg(\tP,e_\tP)=(\tQ,{^\tg e_\tP})\leq (\tD,e_\tD),
$
which, by the uniqueness of Brauer subpairs, is equivalent to showing that $^\tg e_\tP=e_\tQ$.  It is enough to show that $\br_\tQ({^\tg\ell})\leq e_\tQ$, i.e.,
\[
e_\tQ.\br_\tQ({^\tg\ell})\neq 0.
\]
As $(\tQ,e_\tQ)\leq (\tD,e_\tD)$, we already have $\br_\tQ(\ell).e_\tQ=\br_\tQ(\ell)\neq 0$.  We then compute
\[
\br_\tQ(\ell).\cG(\tQ).e_\tQ.br_\tQ({^\tg\ell})=\br_\tQ(\ell).e_\tQ.\cG(\tQ).\br_\tQ({^\tg\ell})=\br_\tQ(\ell).\cG(\tQ).\br_\tQ({^\tg\ell}),
\]
so it suffices to show that the last $\mathbbm{k}$-space is nonzero.  We have
\[
\br_\tQ(\ell).\cG(\tQ).\br_\tQ({^\tg\ell})=\br_\tQ(\ell.\cG^\tQ\cdot \tg\cdot \ell\cdot \tg^{-1})
\]
and, as $\tg\in\cG^\times\cap{^\varphi\cG^\tP}$, it is easy to check that $\cG^\tQ\cdot \tg={^\varphi\cG^\tP}$.  It is also immediate that 
\[
\ell.{^\varphi\cG^\tP}.\ell={^\varphi\!(\ell.\cG.\ell)^\tP}={^\varphi\!\cS^\tP},
\]
so the $\mathbbm{k}$-space we are interested in is actually $\br_\tQ({^\varphi\!\cS^\tP}\cdot \tg^{-1}).$

This last space contains the element $\br_\tQ(x\cdot \tg^{-1})=\br_\varphi(x).\br_{\varphi^{-1}}\!(\tg^{-1})$, which is nonzero as $x$ is an element of the $(\tD,\tD)$-invariant basis $X$  and $\tg$ is a unit in $\cG$.  This completes the proof of the inclusion $\fF_\tD(\cS)\subseteq\fF_\tD(b)$.

 $\bm{\cF_\tD(b)\subseteq\fF_\tD(\cS)}${\bf:}  Suppose now that $\varphi\in\Iso_{\cF_\tD(b)}(\tP,\tQ)$, so that there is some $\tg\in \tG$ with $^\tg(\tP,e_\tP)\leq (\tD,e_\tD)$ and $\varphi=c_\tg|_\tP$.  We make use of the fact that $\cF_\tD(b)$ is a saturated fusion system, and as such Alperin's fusion theorem applies \cite[Theorem A.10]{BLO2}.  In particular, there are $\cF_\tD(b)$-centric subgroups $\tR_1,\ldots,\tR_n\leq \tD$ and automorphisms $\alpha_i\in\Aut_{\cF_\tD(b)}(\tR_i)$ such that $\tP\leq \tR_1$ and
\[
\varphi=\alpha_n\circ\ldots\circ\alpha_1,
\]
with the obvious restrictions omitted from the notation.  Thus there are group elements $\tg_1,\ldots,\tg_n\in \tG$ such that $^{\tg_i}(\tR_i,e_{\tR_i})=(\tR_i,e_{\tR_i})$ and $\alpha_i=c_{\tg_i}|_{\tR_i}$.  As each $\tR_i$ is centric in $\cF_\tD(b)$, there is a unique local point $\widehat \varepsilon_i\in\cL\cP(\cG^{\tR_i})$ such that $(\tR_i)_{\widehat\varepsilon_i}\leq \tD_{[\ell]}$ in $s(\tD,\cG)$ (cf. \cite[Proposition 41.1]{ThevenazBook}), so we must have $^{\tg_i}((\tR_i)_{\widehat{\varepsilon}_i})=((\tR_i)_{\widehat{\varepsilon}_i})$. Equivalently, if we set $\overline\varepsilon_i:=\widehat \varepsilon_i\cap\cB$, then $\overline \varepsilon_i$ is the unique local point of $\cB^{\tR_i}$, and we have $^{\tg_i}((\tR_i)_{\overline\varepsilon_i})=((\tR_i)_{\overline\varepsilon_i})$.  Therefore $\alpha_i\in\Aut_{\fL_\tD^\tG(\cB)}((\tR_i)_{\overline \varepsilon_i})$.

By the uniqueness of $\overline\varepsilon_i$ in $\cL\cP(\cB^{\tR_i})$, it follows that $\varepsilon_i:=\overline\varepsilon_i\cap\cS\neq\emptyset$, so we have $(\tR_i)_{\varepsilon_i}\in s(\tD,\cS)$.  Proposition \ref{prop:local_Puig_categories_equivalent} implies that $\Aut_{\fL_\tD^\tG(\cB)}((\tR_i)_{\overline\varepsilon_i})=\Aut_{\fL_\tD^\times(\cS)}((\tR_i)_{\varepsilon_i})$, and as we have already noted that $\alpha_i$ lies in the former, there must be some unit $u_i\in\cS^\times$ realizing its inclusion in the latter.  Set $u=u_n.\ldots.u_2. u_1$.  Then for any local pointed subgroup $\tP_\gamma\in s(\tD,\cS)$ such that $\tP_\gamma\leq (\tR_1)_{\varepsilon_1}$ (which again exists by the uniqueness of $\varepsilon_1$ in $\cL\cP(\cS^{\tR_1})$), we have $\delta:={^{u}\gamma}\in\cL\cP(\cS^\tQ)$ and $u$ realizes the $\cS^\times$-isofusion $\varphi:\tP_\gamma\xrightarrow\cong\tQ_\delta$. Lemma \ref{lem:Puig_fusion_containment} then yields $\varphi\in\Iso_{\fF_\tD(\cS)}(\tP,\tQ)$, and the result is proved.
\end{proof}
\end{theorem}

In light of our overarching goal, it is worth reformulating the containment $\fF_\tD(\cS)\subseteq\cF_\tD(b)$ of Theorem \ref{thm:source_algebras_are_divisible} as a result in its own right:

\begin{cor}\label{cor:generation_of_source_basis}
If $X$ is a $(\tD,\tD)$-invariant $\cO$-basis of the source algebra $\cS$, then $X$ is $\cF_\tD(b)$-generated as a $(\tD,\tD)$-biset.
\end{cor}

As noted above, Corollary \ref{cor:generation_of_source_basis} is a known property of source algebras, but the reverse inclusion $\cF_\tD(b)\subseteq\fF_\tD(\cS)$ appears to us to be novel.  It generalizes a result of Linckelmann \cite[7.8]{LinckelmannOnSpendidDerived}, which asserts, in different language, that if $\varphi\in\Iso_{\cF_\tD(b)}(\tP,\tQ)$ and $\tQ$ is fully $\cF_\tD(b)$-centralized, then $\varphi\in\Hom_{\fF_\tD(\cS)}(\tP,\tQ)$.  This earlier result implied that a $(\tD,\tD)$-invariant $\cO$-basis of $\cS$ determined $\cF_\tD(b)$ as the fusion system \emph{generated} by morphisms appearing in point-stabilizers; we have now shown that $\cF_\tD(b)$ is recovered in this manner on the nose.

The second comment concerns generalizing Theorem \ref{thm:source_algebras_are_divisible} itself.  Our ability to define a fusion system for the block $b$ relies on having a well-defined notion of containment of Brauer pairs, and our ability to define this containment in terms of local pointed groups relies in part on the minor miracle that if $\ell\in\cB^\tD$ is a source idempotent, then $\br_\tP(\ell)$ is contained in a unique block $e_\tP$ of $\cB(\tP)$ for all $\tP\leq\tD$ (cf. Proposition \ref{prop:Brauer_partial_order_properties}\ref{prop:Brauer_partial_order_properties_ii}).  Linckelmann \cite[\S4]{LinckelmannTrivialSource} takes this as motivation to define an \emph{almost-source idempotent} for the block $b$ to be any $\widehat{\ell}\in\cB^\tD$ such that $\br_\tD(\widehat\ell)\neq 0$ and $\br_\tD(\widehat{\ell})$ is contained in a unique block of $\cB(\tP)$ for all $\tP\leq\tD$; the corresponding \emph{almost-source algebra}  is $\widehat{\cS}:=\widehat{\ell}.\cB.\widehat{\ell}$. Essentially, the primitivity conditition of the source idempotent is relaxed, but not so far that we lose the ability to describe the containment of Brauer pairs in terms of the natural partial order on the pointed Brown poset.  We then have:

\begin{cor}
Any almost-source algebra $\widehat{\cS}$ of $b$ satisfies $\fF_\tD(\widehat\cS)=\cF_\tD(b)$.  In particular, almost-source algebras are divisible $\tD$-algebras.
\begin{proof}
  To show $\fF_\tD(\widehat{\cS})\subseteq\cF_\tD(b)$, the proof of Theorem \ref{thm:source_algebras_are_divisible} applies exactly with the almost-source idempotent $\widehat\ell$ in place of the source idempotent ${\ell}$, as the same uniqueness properties regarding Brauer pairs are satisfied by both.  Conversely, to see that $\cF_\tD(b)\subseteq\fF_\tD(\widehat{\cS})$, note that the assumption that $\br_\tD(\widehat{\ell})\neq 0$ implies that there is some source idempotent $\ell\leq\widehat{\ell}$, so the almost-source algebra $\widehat{\cS}$ contains a source algebra $\cS$ as a corner subalgebra.  Then $\fF_\tD(\cS)\subseteq\fF_\tD(\widehat{\cS})$, which combined with the containment $\cF_\tD(b)\subseteq\fF_\tD(\cS)$ from before gives our desired result.
\end{proof}
\end{cor}

We close this section with the observation that, in the presence of a unital basis, all isomorphisms of the fixed-point fusion presystem $\fF_\tS(\cA)$ can be viewed as $\cA^\times$-isofusions in all possible ways.
 
 \begin{prop}\label{prop:unital_fusion_Puig_lifting}
Let $\cA$ be a divisible $\tS$-algebra with unital $(\tS,\tS)$-invariant $\cO$-basis $Y$.  Then for every isomorphism $\varphi\in\Iso_{\fF_\tS(\cA)}(\tP,\tQ)$ and every local point $\gamma\in\cL\cP(\cA^\tP)$,  there is a unique local point $\delta\in\cL\cP(\cA^\tQ)$ such that $\varphi\in\Iso_{\fL_\tS^\times(\cA)}(\tP_\gamma,\tQ_\delta)$.  

In particular, for each $\gamma\in\cL\cP(\cA^\tP)$ and $\delta\in\cL\cP(\cA^\tQ)$, we have
\[
\Iso_{\fF_\tS(\cA)}(\tP,\tQ)=\coprod_{\delta'\in\cL\cP(\cA^\tQ)}\Iso_{\fL_\tS^\times(\cA)}(\tP_\gamma,\tQ_{\delta'})=\coprod_{\gamma'\in\cL\cP(\cA^\tP)} \Iso_{\fL_\tS^\times(\cA)}(\tP_{\gamma'},\tQ_\delta).
\]
\begin{proof}
Choose  $y\in Y$ realizing $\varphi\in\Iso_{\fF_\tS(\cA)}(\tP,\tQ)$, so that by definition $y\in{^\varphi\!\!\cA^\tP}$. As $y\in\cA^\times$, we also have $y^{-1}\in{^{\varphi^{-1}}\!\!\!\cA^\tQ}$.  Fix some $i_\tP\in\gamma$ and set $j_\tQ:={^yi_\tP}$. It is immediate that $j_\tQ$ is a primitive local idempotent of $\cA^\tQ$; set  $\delta:=[j_\tQ]\in\cL\cP(\cA^\tQ)$.  Then
\[
\qquad s_\varphi:=j_\tQ.y.i_\tP\in j_\tQ.{^\varphi\!\!\cA^\tP}.i_\tP\qquad\textrm{and}\qquad
t_{\varphi^{-1}}:=i_\tP.y^{-1}.j_\tQ\in i_\tP.{^{\varphi^{-1}}\!\!\!\cA^\tQ}.j_\tQ
\]
satisfy the conditions of  Lemma \ref{lem:Puig_morphisms_via_transpotents}, and thus we have $\varphi\in\Iso_{\fL^\times(\cA)}(\tP_\gamma,\tQ_\delta)$.  The uniqueness of $\delta$ is then a direct application of Lemma \ref{lem:Puig_category_unique_target}.

The first equality of the last claim follows immediately from the above.  The second equality is a consequence of the divisibility of the categories $\fF_\tS(\cA)$ and $\fL_\tS^\times(\cA)$, which allows us to run the same argument after taking the inverse of each element in $\Iso_{\fF_\tS(\cA)}(\tP,\tQ)$ to obtain $\Iso_{\fF_\tS(\cA)}(\tQ,\tP)$.
\end{proof}
\end{prop}

\section{Unital bases}\label{sec:unital_bases}

Let $\tS$ be a $p$-group and $\cA$ be a bifree $\tS$-algebra.   In this section we investigate the implications of the existence of a unital $(\tS,\tS)$-invariant $\cO$-basis $Y\subseteq\cA^\times$.  Note that the importance of such a basis has already been established in the previous discussion (e.g., Lemma \ref{lem:divisibility_sufficiency} and Propositions \ref{prop:local_Puig_morphisms_realized_by_unital basis} and \ref{prop:unital_fusion_Puig_lifting}).

We begin by examining the effect  a unit  lying in a particular fixed-point submodule has on twisting the biset structure of $Y$.

\begin{lemma}\label{lem:units_in_fixed_point_submodules}
Let $\cA$ be a bifree $\tS$-algebra with $(\tS,\tS)$-invariant $\cO$-basis $Y$.  If $\tP\leq\tS$ is a subgroup and $\varphi:\tP\to\tS$ is such that $\cA^\times\cap{^\varphi\!\!\cA^\tP}\neq\emptyset$, then ${^\varphi_\tP Y_\tS}\cong{_\tP Y_\tS}$ as $(\tP,\tS)$-bisets.
\begin{proof}
Pick $u\in\cA^\times\cap{^\varphi\!\!\cA^\tP}$  and set $Y':=u^{-1}.Y$.  It is immediate that $u^{-1}\in\cA^\times\cap{^{\varphi^{-1}}\!\!\!\cA^{\varphi\tP}}$, so for any  $\tp\in\tP$, $y\in Y$, and $\ts\in\tS$, we have
\[
\tp\cdot (u^{-1}.y)\cdot\ts=u^{-1}\cdot\varphi(\tp)\cdot y\cdot\ts\in u^{-1}.Y,
\]
and thus $Y'$  is a $(\tP,\tS)$-biset.  Moreover, as $Y'$ is formed by multiplying an $\cO$-basis of $\cA$ by a unit,  it is also an $\cO$-basis of $\cA$.   Corollary \ref{cor:basis_well_defined}\ref{cor:basis_well_defined_i} then implies that  $_\tP Y'_\tS\cong {_\tP Y_\tS}$ as $(\tP,\tS)$-bisets.

Now consider the map $f:Y'\to Y:y'\mapsto u.y'$.  For $\tp\in \tP$, $y'\in Y'$,  and $\ts\in \tS$ we have
\[
f(\tp\cdot y'\cdot \ts)=u\cdot \tp\cdot y'\cdot \ts=\varphi(\tp)\cdot u.y'\cdot \ts=\varphi(\tp)\cdot f(y')\cdot \ts,
\]
showing that $f$ is an isomorphism $Y'\cong{^\varphi_\tP Y_\tS}$ of $(\tP,\tS)$-bisets.  Composing these isomorphisms yields ${^\varphi_\tP Y_\tS}\cong{_\tP Y_\tS}$, as desired.
\end{proof}
\end{lemma}

When $\cA$ is a divisible $\tS$-algebra, the existence of a unital $(\tS,\tS)$-invariant basis has the following immediate consequence.

\begin{prop}\label{prop:basis_of_units_implies_conj1}
If $\cA$ is a divisible $\tS$-algebra with unital $(\tS,\tS)$-invariant $\cO$-basis $Y$, then $Y$ is $\fF_\tS(\cA)$-stable.
\begin{proof}
Recall that the divisibility of $\cA$ implies that the fixed-point fusion presystem is actually a fusion system on $\tS$, so it makes sense to speak of $\fF_\tS(\cA)$-stability.

By the definition of $\fF_\tS(\cA)$, for any $\varphi\in\Hom_{\fF_\tS(\cA)}(P,S)$ there is some $y\in{^\varphi Y^P}$.  As $Y\subseteq\cA^\times$ consists of units, the hypotheses of Lemma \ref{lem:units_in_fixed_point_submodules} apply, so we have $_\tP^\varphi Y_\tS\cong{_\tp Y_\tS}$ as $(\tP,\tS)$-bisets, which is exactly the definition of $\fF_\tS(\cA)$-stability.
\end{proof}
\end{prop}

In particular, we have reproved in slightly greater generality a result known to Linckelmann and Webb, which can be found in \cite[Proposition 8.7.11]{LinckelmannBookII}:

\begin{cor}\label{cor:strong_1_implies_main}
Conjecture \ref{conj:main_strong_i} implies Conjecture \ref{conj:main}.
\begin{proof}
By the discussion in \S\ref{subsec:state_of_the_art}, we must only verify the $\cF_\tD(b)$-stability of the $(\tD,\tD)$-invariant $\cO$-basis $X$ of $\cS$.  If $X$ can be chosen to be unital, Proposition \ref{prop:basis_of_units_implies_conj1} implies that $X$ is $\fF_\tD(\cS)$-stable, but $\fF_\tD(\cS)=\cF_\tD(b)$ by Theorem \ref{thm:source_algebras_are_divisible}.  The result follows.
\end{proof}
\end{cor}

The other structural conditions on $\cA$ we will introduce are, at  heart, different sets of hypotheses that  guarantee the existence of a unital $(\tS,\tS)$-invariant $\cO$-basis.  The most basic of these asks that all relevant twisted-diagonal fixed-point submodules of $\cA$ contain a unit. We close this section with an exploration of this property, beginning with some basic linear algebra and facts  relating the $\cO$-algebra $\cA$ to the residue $\kk$-algebra $\overline\cA$.

\begin{lemma}\label{lem:lin_alg_basis}
Let $V$ be a $\kk$-vector space with $\kk$-basis $Z$.  For $z\in Z$, $v\in V$,  and $k\in\kk$,  let $Z_k$ be $Z$ with $z$ replaced by $z+kv$.  Then $Z_k$ is a $\kk$-basis of $V$ for all but finitely many $k$.
\begin{proof}
Let $M_k$ be the $(\dim_\kk V\times\dim_\kk V)$-matrix whose columns represent the elements of $Z_k$ in terms of the basis $Z$.  Clearly each column has a unique nonzero entry, which is $1$, except for the column corresponding to $z+kv$.  It is also clear that $\det(M_k)$ is polynomial in $k$, and is not identically 0 as evaluation at $k=0$  shows.  Thus $\det(M_k)\neq 0$ for all but finitely many $k\in\kk$,  and for any such we have that $Z_k$ is a $\kk$-basis of $V$.
\end{proof}
\end{lemma}

\begin{lemma}\label{lem:k_to_O_lifting}
Let $Y\subseteq\cA$ be a linearly independent set whose image $\overline Y\subseteq\overline\cA$ is also linearly independent.  Then $Y$ is an $\cO$-basis for $\cA$ if and only if $\overline Y$ is a $\kk$-basis for $\overline\cA$.
\begin{proof}
Clearly if $Y$ is an $\cO$-basis then $\overline Y$ is a $\kk$-basis, so we assume the second condition.  Let $Z\subseteq\cA$ be an $\cO$-basis of $\cA$, and let $M$ be the $(\textrm{rank}_\cO\phantom{.}\cA\times\textrm{rank}_\cO\phantom{.}\cA)$-matrix whose columns represent the elements of $Y$ in terms of $Z$.  As $\overline Z$ and $\overline Y$ are both $\kk$-bases for $\overline\cA$, the residue matrix $\overline M$ represents a change of basis operation, and thus is invertible.  In particular, $0\neq\det(\overline M)=\overline{\det(M)}$, so $\det(M)\in\cO^\times$ and $M$  represents an invertible linear transformation.  It follows that $Y$ is a $\cO$-basis for $\cA$.
\end{proof}
\end{lemma}

\begin{lemma}\label{lem:lin_alg_units}
For  $a\in\cA$, $u\in\cA^\times$, and $\lambda\in\cO$, the expression $a+\lambda u$ is a unit in $\cA^\times$ for all but finitely many values of $\overline\lambda\in\kk$.
\begin{proof}
Let $\ell_-:\cA\to\End_\cO(\cA)$ and $\overline\ell_-:\overline\cA\to\End_\kk(\overline\cA)$ be the left regular representations.  We have $a+\lambda u\in\cA^\times$ if and only if $\det(\ell_{a+\lambda u})\in\cO^\times$ if and only if $\det(\overline\ell_{\overline{a+\lambda u}})\neq 0\in\kk$ if and only if $\overline a+\overline \lambda\overline u\in\overline \cA^\times$.  For $\overline\lambda\neq 0$ we clearly have $\overline a+\overline\lambda\overline u\in\overline A^\times$ if and only if $\overline\lambda^{-1}\overline a+\overline u\in\cA^\times$, so we must show that  $\det(\overline\ell_{k\overline a+\overline u})\neq 0$ for all by finitely many $k\in \kk$.  Clearly this determinant is  polynomial in $k$ and is not identically $0$ as evaluation at $k=0$ and the fact that $\overline u$ is a unit shows.  The result follows.
\end{proof}
\end{lemma}

\begin{prop}\label{prop:unit_stabilizers_implies_unit_basis}
Let $\cA$ be a bifree $\tS$-algebra.  Then $\cA$ possesses a unital $(\tS,\tS)$-invariant $\cO$-basis if and only if for all $\varphi\in\Hom_{\fF_\tS(\cA)}(\tP,\tS)$ we have $\cA^\times\cap{^\varphi\!\!\cA^\tP}\neq\emptyset$.
\begin{proof}
Let $Y$ be an  $(\tS,\tS)$-invariant $\cO$-basis of $\cA$.  If $Y\subseteq\cA^\times$, and $\varphi\in\Hom_{\fF_\tS(\cA)}(\tP,\tQ)$ is realized by $y\in Y$, then $y\in\cA^\times\cap{^\varphi\!\!\cA^\tP}$, so the ``only if'' direction is immediate.

We now prove the ``if'' implication. Consider  $y\in Y$: If $y\in\cA^\times$ is a unit, then so is the entire $(\tS,\tS)$-orbit $\tS\cdot y\cdot\tS$, in which case we do nothing.  If $y\notin\cA^\times$, we will replace the $(\tS,\tS)$-orbit of $y$ by an $(\tS,\tS)$-orbit of units, in a manner still yielding an $\cO$-basis for $\cA$.

Write $\Stab(y)=\Delta(\varphi,\tP)$, so by hypothesis there is some  $u\in\cA^\times\cap{^\varphi\!\!\cA^\tp}$.  In particular, we have $\Stab(y)\leq\Stab(u)$.  It follows that for all $\lambda\in\cO$, the map of $(\tS,\tS)$-orbits
\[
\tS\cdot y\cdot\tS\to\tS\cdot(y+\lambda u)\cdot\tS:\ts_1\cdot y\cdot \ts_2 \mapsto\ts_1\cdot (y+\lambda u)\cdot\ts_2
\]
is well-defined, as $\ts_1\cdot y\cdot \ts_2=\ts_1'\cdot y\cdot \ts_2'$ implies $((\ts_1')^{-1}\ts_1,\ts_2'\ts_2^{-1})\in\Stab(y)\leq\Stab(u)$ so that $\ts_1\cdot(y+\lambda u)\cdot\ts_2=\ts_1'\cdot(y+\lambda u)\cdot\ts_2'$.  

For a fixed $y'\in\tS\cdot y\cdot \tS$, let $Y'_\lambda$ be $Y$ with $y'$ replaced by $y'+\lambda u$.  In the residue algebra $\overline\cA$, Lemma \ref{lem:lin_alg_basis} tells us that $\overline {Y'_\lambda}$ is a $\kk$-basis for all but finitely many $\overline\lambda\in\kk$.  For any  $\overline\lambda$ that does yield a $\kk$-basis,  $Y'_\lambda$ is an $\cO$-basis of $\cA$ by Lemma \ref{lem:k_to_O_lifting}.  Letting $y'$ range over the orbit of $y$, and setting $Y_\lambda$ to be $Y$ with the entire $(\tS,\tS)$-orbit $\tS\cdot y\cdot\tS$ replaced by $\tS\cdot(y+\lambda u)\cdot\tS$, we have shown that $Y_\lambda$ is an $\cO$-basis of $\cA$ for all but finitely many values of $\overline\lambda\in\kk$.

On the other hand, Lemma \ref{lem:lin_alg_units} gives $Y_\lambda\subseteq\cA^\times$ for all but finitely many $\overline\lambda\in\kk$.  The infinitude of the residue field implies that there is some $\lambda\in\cO$ such that $Y_\lambda$ is  an $(\tS,\tS)$-invariant $\cO$-basis of $\cA$ with one more $(\tS,\tS)$-orbit of units than possessed by $Y$.  Proceeding orbit-by-orbit yields the desired unital $(\tS,\tS)$-invariant basis, and we are done.
\end{proof}
\end{prop}

In what follows, Proposition \ref{prop:unit_stabilizers_implies_unit_basis} will be our main point of connection between the existence of a unital $(\tS,\tS)$-invariant $\cO$-basis of $\cA$ and the other properties we consider.

\section{Twisted units}\label{sec:twisted_units}

In this section, we investigate the behavior of a unital invariant basis on the level of twisted Brauer quotients.  We  characterizatize the existence of such a basis in these terms, which proves the equivalence of Conditions  \ref{thm:main_i} and \ref{thm:main_ii} of Theorem \ref{thm:main}.

Suppose that $\cA$ is a divisible $\tS$-algebra with unital $(\tS,\tS)$-invariant $\cO$-basis $Y$.  If $y\in Y$ has point-stabilizer $\Delta(\varphi,\tP)$,   Lemma \ref{lem:Brauer_quotient_dimensions} tells us that $\br_\varphi(y)$ is nonzero in $\cA(\varphi)$.  On the other hand, the proof of Lemma \ref{lem:divisibility_sufficiency} implies that, as $y^{-1}$ is a unit, we have  $\br_{\varphi^{-1}}(y^{-1})$ is nonzero in $\cA(\varphi^{-1})$.  We can use these facts to compare different twisted Brauer quotient modules.  If we set $\tQ:=\varphi(\tP)$, then
\[
\cA(\tP)\to\cA(\varphi):a_\tP\mapsto \br_\varphi(y).a_\tP\qquad\textrm{and}\qquad
\cA(\tQ)\to\cA(\varphi):a_\tQ\mapsto a_\tQ.\br_\varphi(y)
\]
are $\mathbbm{k}$-linear isomorphisms, with inverses given by
\[
\cA(\varphi)\to\cA(\tP):a_\varphi\mapsto \br_{\varphi^{-1}}(y^{-1}).a_\varphi\qquad\textrm{and}\qquad
\cA(\varphi)\to\cA(\tQ):a_\varphi\mapsto a_\varphi.\br_{\varphi^{-1}}(y^{-1}),
\]
respectively.  Moreover, the maps involving $\cA(\tP)$ are actually morphisms of right $\cA(\tP)$-modules, while those involving $\cA(\tQ)$ are left $\cA(\tQ)$-module maps.  In other words,  when we view $\cA(\varphi)$ as an $(\cA(\tQ),\cA(\tP))$-bimodule, the unital basis element $y$ gives rise to a simultaneous generator $\br_\varphi(y)$ for $\cA(\varphi)$ as a regular left $\cA(\tQ)$-module and regular right $\cA(\tP)$-module.  In particular, the $\kk$-dimensions of $\cA(\tP)$, $\cA(\tQ)$, and $\cA(\varphi)$ are all equal. It is not hard to see that this line of reasoning leads to the conclusion that $Y$ is $\fF_\tS(\cA)$-stable, giving us a second proof of Proposition \ref{prop:basis_of_units_implies_conj1}.  

We now set out to investigate the converse problem:  If, for all $\varphi\in\Iso_{\fF_\tS(\cA)}(\tP,\tQ)$, the $(\cA(\tQ),\cA(\tP))$-bimodule $\cA(\varphi)$ is simultaneously biregular  in the sense above, does it follow that there exists an $(\tS,\tS)$-invariant $\cO$-basis of units of $\cA$?

We outline the basic flow of the argument.  In Definition \ref{def:twisted_units} below, we introduce the notion of a ``twisted unit,'' which is designed to mimic the image of a unit of $\cA$ in the twisted Brauer quotient $\cA(\varphi)$.   If we then assume that for all $\varphi\in\Iso_{\fF_\tS(\cA)}(\tP,\tQ)$ the set of twisted units in $\cA(\varphi)$ is nonempty, we will show that every such twisted unit can be lifted to a unit in $\cA$ whose point-stabilizer contains $\Delta(\varphi,\tP)$.  Once this is established, Proposition \ref{prop:unit_stabilizers_implies_unit_basis} will then give us the existence of a unital $(\tS,\tS$)-invariant $\cO$-basis.

The lifting process from twisted units to units is a bit involved, but in brief:
\begin{enumerate}
\item[(1)] Fix $\varphi\in\Iso_{\fF_\tS(\cA)}(\tP,\tQ)$ and primitive idempotent decompositions
\[
1_{\cA^\tP}=\sum_{i\in I}i\qquad\textrm{and}\qquad 1_{\cA^\tQ}=\sum_{j\in J}j
\]
in the fixed-point subalgebras $\cA^\tP$ and $\cA^\tQ$, respectively.
\item[(2)] Establish a bijection $\sigma:I\to J$ with the property that for all $i\in I$, there are elements
\[
u_i\in\sigma(i)_.{^\varphi\!\!\cA^\tP}.i\qquad\textrm{and}\qquad v_i\in i.{^{\varphi^{-1}}\!\!\!\cA^\tQ}.\sigma(i)
\]
such that $v_i.u_i=i$ and $u_i.v_i=\sigma(i)$.  Note that $u_i.v_{i'}=0=v_{i'}.u_i$ for $i\neq i'$ as the idempotents of $I$ and $J$ are mutually orthogonal.
\item[(3)] Set
\[
u=\sum_{i\in I}u_i\qquad\textrm{and}\qquad v=\sum_{i\in I} v_j.
\]
We have $v.u=\sum\limits_{i\in I}i=1_{\cA^\tP}=1_\cA$ and $u.v=\sum\limits_{j\in J}j=1_{\cA^\tQ}=1_\cA$, so that $u$ and $v$ are mutually inverse units of $\cA$.  As $u\in {^\varphi\!\!\cA^\tP}$ by construction, the goal follows.
\end{enumerate}

With our plan of attack established, we now give our main definition. 

\begin{definition}\label{def:twisted_units}
Let $\cA$ be a divisible $\tS$-algebra and $\varphi\in\Iso_{\fF_\tS(\cA)}(\tP,\tQ)$. A \emph{twisted unit} of $\varphi$ is any element $u_\varphi\in\cA(\varphi)$ such that there exists  $u_\varphi^\dag\in\cA(\varphi^{-1})$ with $u_\varphi^\dag.u_\varphi=1_{\cA(\tP)}$. The element $u_\varphi^\dag$ is a \emph{twisted inverse} of $u_\varphi$.  The set of all twisted units of $\varphi$ is denoted $\cA(\varphi)^\times$.
\end{definition}

We say that $\cA$ \emph{has all twisted units} if each isomorphism of $\fF_\tS(\cA)$ has a twisted unit.

Before investigating the properties of twisted units in more detail, let us make note of an equivalent formulation of $\cA$'s having all twisted units.

\begin{prop}\label{prop:PRS_alg_equiv}
The divisible $\tS$-algebra $\cA$ has all twisted units if and only if, for all composable isomorphisms $\tP\xrightarrow\varphi \tQ\xrightarrow\psi \tR$ of $\fF_\tS(\cA)$, the bilinear pairing
\[
\cA(\psi)\times\cA(\varphi)\to\cA(\psi\varphi)
\]
is a surjection.
\begin{proof}
Suppose that $\cA$ has all twisted units and we are given a composable pair of isomorphisms $\tP\xrightarrow\varphi \tQ\xrightarrow\psi \tR$ as well as some $a_{\psi\varphi}\in\cA(\psi\varphi)$.  As $\cA(\varphi)^\times\neq\emptyset$, we may choose $u_\varphi\in\cA(\varphi)$ and $u_\varphi^\dag\in\cA(\varphi^{-1})$ such that $u_\varphi^\dag.u_\varphi=1_{\cA(\tP)}$.  We then have
\[
a_{\psi\varphi}=a_{\psi\varphi}.1_{\cA(\tP)}=(a_{\psi\varphi}.u_\varphi^\dag).u_\varphi,
\]
which factors $a_{\psi\varphi}$ as a product of elements from $\cA(\psi)$ and $\cA(\varphi)$, as desired.

The converse is proved by taking $\psi=\varphi^{-1}$ and choosing $a_\varphi\in\cA(\varphi)$, $a_{\varphi^{-1}}\in\cA(\varphi^{-1})$ that satisfy $a_{\varphi^{-1}}.a_\varphi=1_{\cA(\varphi)}$.
\end{proof}
\end{prop}

We here collect various facts on twisted units that will be essential in the sequel.  Note that the overarching theme is that twisted units behave very much as if they were actual units, both in intrinsic terms related to their composition, uniqueness of inverses, etc., as well as their extrinsic relations to Brauer quotients.  In the following we fix a pair of composable isomorphisms $\tP\xrightarrow\varphi \tQ\xrightarrow\psi \tR$ in $\fF_\tS(\cA)$.

We first note that twisted units are closed under multiplication, where defined:

\begin{lemma}\label{lem:twisted_units_properties_i}
If $u_\varphi\in\cA(\varphi)^\times$ and $u_\psi\in\cA(\psi)^\times$, then $u_\psi.u_\varphi\in\cA(\psi\varphi)^\times$.
\begin{proof}
If $u_\varphi^\dag\in\cA(\varphi^{-1})$ and $u_\psi^\dag\in\cA(\psi^{-1})$ are twisted inverses of $u_\varphi$ and $u_\psi$, we have
\[
(u_\varphi^\dag.u_\psi^\dag).(u_\psi.u_\varphi)=u_\varphi^\dag.1_{\cA(Q)}.u_\varphi=1_{\cA(P)},
\]
so that $u_\varphi^\dag.u_\psi^\dag$ is a twisted inverse of $u_\psi.u_\varphi$.
\end{proof}
\end{lemma}

As in the introduction to this section, multiplication by a twisted unit induces an isomorphism between twisted Brauer quotients:

\begin{lemma}\label{lem:twisted_units_properties_ii}
If $u_\varphi\in\cA(\varphi)^\times$ has twisted inverse $u_\varphi^\dag\in\cA(\varphi^{-1})$, then the maps
\[
\cA(\tP)\to \cA(\varphi):a_\tP\mapsto u_\varphi.a_\tP\qquad\textrm{and}\qquad
\cA(\varphi)\to \cA(\tQ):a_\varphi\mapsto a_\varphi.u_\varphi^\dag
\]
are isomorphisms of $\mathbbm{k}$-vector spaces.  In particular, $u_\varphi$ is a simultaneous generator of $\varphi(\cA)$ as a regular right $\cA(\tP)$-module and regular left $\cA(\tQ)$-module.
\begin{proof}
The given maps both have left inverses, given by left multiplication by $u_\varphi^\dag$ and right mutliplication by $u_\varphi$, respectively.  Thus each is a $\mathbbm{k}$-linear injection, and
\[
\dim_\kk\cA(\tP)\leq\dim_\kk\cA(\varphi)\leq\dim_\kk\cA(\tQ).
\]
These dimension bounds depend only on the existence of a twisted unit in $\cA(\varphi)^\times$.  As $\cA$ has all twisted units, we have $\cA(\varphi^{-1})^\times\neq\emptyset$, so the same argument yields
\[
\dim_\kk\cA(\tQ)\leq\dim_\kk\cA(\varphi^{-1})\leq\dim_\kk\cA(\tP).
\]
Thus  the $\mathbbm{k}$-dimensions of $\cA(\tP)$, $\cA(\tQ)$, and $\cA(\varphi)$ are all equal, so the given maps are both $\kk$-linear isomorphisms.  The final statement is immediate.
\end{proof}
\end{lemma}

Twisted inverses are twisted units in their own right, and are unique:

\begin{lemma}\label{lem:twisted_units_properties_iii}
If $u_\varphi\in\cA(\varphi)^\times$ has twisted inverse $u_\varphi^\dag\in\cA(\varphi^{-1})$, then $u_\varphi.u_\varphi^\dag=1_{\cA(\tQ)}$, so that $u_\varphi^\dag\in\cA(\varphi^{-1})^\times$ and $u_\varphi$ is a twisted inverse of $u_\varphi^\dag$.  Moreover, $u_\varphi^\dag$ is the unique element of $\cA(\varphi^{-1})$ such that $u_\varphi^\dag.u_\varphi=1_{\cA(\tP)}$.
\begin{proof}
By Lemma \ref{lem:twisted_units_properties_ii}, there is some $a_\varphi\in\cA(\varphi)$ such that $a_\varphi.u_\varphi^\dag=1_{\cA(\tQ)}$.  Right multiplying by $u_\varphi$ yields $a_\varphi=u_\varphi$, proving the first claim.  

Now suppose that $v_{\varphi^{-1}}\in\cA(\varphi^{-1})$ also satisfies $v_{\varphi^{-1}}.u_\varphi=1_{\cA(\tP)}$, so $v_{\varphi^{-1}}.u_\varphi.u_\varphi^\dag=u_\varphi^\dag$.  As we've just verified that $u_\varphi.u_\varphi^\dag=1_{\cA(\tQ)}$, we have $v_{\varphi^{-1}}=u_\varphi^\dag$, as desired.
\end{proof}
\end{lemma}

A consequence of Lemma \ref{lem:twisted_units_properties_iii} is that it makes sense to speak of \emph{the} twisted inverse of a twisted unit, and we will do so from now on.

For $\varphi\in\Iso_{\fF_\tS(\cA)}(\tP,\tQ)$, the set $\cA(\varphi)^\times$ of twisted units of $\varphi$  is naturally a $(\cA(\tQ)^\times,\cA(\tP)^\times)$-biset.  In fact, this biset is biregular:

\begin{lemma}\label{lem:twisted_units_properties_iv}
If $u_\varphi,v_\varphi\in\cA(\varphi)^\times$ are  twisted units, then there exist unique $x_\tP\in\cA(\tP)^\times$ and $x_\tQ\in\cA(\tQ)^\times$ such that $u_\varphi.x_\tP=v_\varphi=x_\tQ.u_\varphi$.
\begin{proof}
 Consider  $u_\varphi^\dag\in\cA(\varphi^{-1})^\times$.  By Lemma \ref{lem:twisted_units_properties_i}, both the products $x_\tP:=u_\varphi^\dag. v_\varphi\in\cA(\tP)$ and $x_\tQ:=v_\varphi.u_\varphi^\dag\in\cA(\tQ)$ are twisted units.  Both $\cA(\tP)$ and $\cA(\tQ)$ are actual algebras, where the notions of unit and twisted unit coincide.  It is then immediate that our choices of $x_\tP$ and $x_\tQ$ satisfy the desired condition, and uniqueness follows from the fact that the action of a group of units in an algebra on the nonzero elements of a regular module is free.
\end{proof}
\end{lemma}

Finally, `conjugation' by a twisted unit induces an algebra isomorphism between (untwisted) Brauer quotients:

\begin{lemma}\label{lem:twisted_units_properties_v}
If $u_\varphi\in\cA(\varphi)^\times$, then the map
\[
\cA(\tP)\to\cA(\tQ):a_\tP\mapsto u_\varphi.a_\tP.u_\varphi^\dag
\]
is an algebra isomorphism.
\begin{proof}
The inverse map is given by $\cA(\tQ)\to\cA(\tP):a_\tQ\mapsto u_\varphi^\dag.a_\tQ.u_\varphi$, applying Lemma \ref{lem:twisted_units_properties_iii}.  That both maps are algebra morphisms is obvious.
\end{proof}
\end{lemma}

We now begin the process of lifting twisted units to units in $\cA^\times$. Lemma \ref{lem:twisted_units_properties_v} tells us that, if $\cA$ has all twisted units, then for any isomorphism $\varphi\in\Iso_{\fF_\tS(\cA)}(\tP,\tQ)$  the untwisted Brauer quotients $\cA(\tP)$ and $\cA(\tQ)$ are abstractly isomorphic as algebras, and that one such isomorphism is determined by each element of $\cA(\varphi)^\times$.  In particular, a choice of twisted unit determines a bijection between any algebraic structure of $\cA(\tP)$ and the corresponding structure of $\cA(\tQ)$---for instance the primitive idempotents of each.  We would like for these bijections to depend only on the $\fF_\tS(\cA)$-isomorphism $\varphi$ and not the particular choice of twisted unit, but clearly this cannot be done on the level of the primitive idempotents themselves.  However, Lemma \ref{lem:twisted_units_properties_iv} implies that we do obtain a well-defined bijection between the points of these  algebras
\[
\overline\theta_\varphi:\cP(\cA(\tP))\to\cP(\cA(\tQ)):[e]\mapsto[u_\varphi.e.u_\varphi^\dag],
\]
which is independent of the choice of twisted unit $u_\varphi\in\cA(\varphi)^\times$.

As the Brauer map $\br_\tP:\cA^\tP\to\cA(\tP)$ induces a bijection $\cL
\cP(\cA^\tP)\xrightarrow\cong\cP(\cA(\tP))$, we may lift $\overline\theta_\varphi$ to a bijection of local points $\theta_\varphi:\cL\cP(\cA^\tP)\to\cL\cP(\cA^\tQ)$, which is characterizied by the commuting diagram:
\[
\xymatrix{
\cL\cP(\cA^\tP)\ar[r]^-{\theta_\varphi}\ar[d]_{\br_\tP}&\cL\cP(\cA^\tQ)\ar[d]^{\br_\tQ}\\
\cP(\cA(\tP))\ar[r]_-{\overline\theta_\varphi}&\cP(\cA(\tQ))
}
\]
In other words, 
\[
\br_\tQ(\theta_\varphi([i]))=\overline\theta_\varphi(\br_\tP([i]))
\]
for any primitive local idempotent $i\in\cL\cP(\cA^P)$.  

The bijection $\theta_\varphi$ will be our main tool for comparing primitive idempotent decompositions of Brauer quotient algebras.  We will often make use of its alternative charaterization:

\begin{prop}\label{prop:PRS_algebra_transpotent_lifting}
For any $\varphi\in\Iso_{\fF_\tS(\cA)}(\tP,\tQ)$, $\gamma\in\cL\cP(\cA^\tP)$, and $\delta\in\cL\cP(\cA^\tQ)$, we have $\theta_\varphi(\gamma)=\delta$ if and only if for some, and hence all, $i_\tP\in\gamma$ and $j_\tQ\in \delta$, there exist 
\[
s_\varphi\in j_\tQ.{^\varphi\!\!\cA^\tP}.i_\tP\qquad\textrm{and}\qquad t_{\varphi^{-1}}\in i_\tP.{^{\varphi^{-1}}\!\!\!\cA^\tQ}.j_\tQ
\]
 such that $i_\tP=t_{\varphi^{-1}}.s_\varphi$ and $j_\tQ=s_\varphi.t_{\varphi^{-1}}$.
\begin{proof}
Suppose that $s_\varphi$ and $t_{\varphi^{-1}}$ exist and pick some $u_\varphi\in\cA(\varphi)^\times$.  Using bars to denote the image under the appropriate Brauer map, we wish to show that $u_\varphi.\overline i_\tP.u_\varphi^\dag$ and $\overline j_\tQ$ are conjugate in $\cA(\tQ)$.  We  expand
\[
u_\varphi.\overline i_\tP.u_\varphi^\dag=(u_\varphi.\overline t_{\varphi^{-1}}).(\overline s_\varphi.u_\varphi^\dag)\qquad\textrm{and}\qquad 
\overline j_\tQ=\overline s_\varphi.\overline t_{\varphi^{-1}}=
(\overline s_\varphi.u_\varphi^\dag).(u_\varphi.\overline t_{\varphi^{-1}}).
\]
As both $u_\varphi.\overline t_{\varphi^{-1}}$ and $\overline s_\varphi.u_\varphi^\dag$ lie in $\cA(\tQ)$, this shows that $u_\varphi.\overline i_\tP.u_\varphi^\dag$ and $\overline j_\tQ$ are associate, and hence conjugate, idempotents of $\cA(\tQ)$.  The ``only if'' implication is proved.

Conversely, suppose that $\theta_\varphi(\gamma)=\delta$, so that $[u_\varphi.\overline i_\tP.u_\varphi^\dag]=[\overline j_\tQ]$ for  $u_\varphi\in\cA(\varphi)^\times$.  By Lemma \ref{lem:twisted_units_properties_iv}, we may choose the twisted unit  so that   $u_\varphi.\overline i_\tP.u_\varphi^\dag=\overline j_\tQ$ on the nose.  Let $\widehat u_\varphi\in{^\varphi\!\!\cA^\tP}$ and $\widehat{u}_\varphi^\dag\in{^{\varphi^{-1}}\!\!\!\cA^\tQ}$ be lifts $u_\varphi$ and $u_\varphi^\dag$, respectively.  Consider the elements 
\[
s_\varphi:=j_\tQ.\widehat u_\varphi.i_\tP\in{^\varphi\!\!\cA^\tP}\qquad\textrm{and}\qquad
\widetilde t_{\varphi^{-1}}:=i_\tP.\widehat u_\varphi^\dag.j_\tQ\in{^{\varphi^{-1}}\!\!\!\cA^\tQ},
\]
and observe that
\[
\br_\tQ(s_\varphi.\widetilde t_{\varphi^{-1}})=\overline j_\tQ.u_\varphi.(\overline i_\tP)^2.u_\varphi^\dag.\overline j_\tQ=(\overline j_\tQ)^4=\overline j_\tQ,
\]
so that $s_\varphi.\widetilde t_{\varphi^{-1}}$ differs from $j_\tQ$ by an element $x\in\ker(\br_\tQ)$: $j_\tQ=s_\varphi.\widetilde t_{\varphi^{-1}}+x$.  As both $j_\tQ$ and $s_\varphi.\widetilde t_{\varphi^{-1}}$ lie in the local algebra $j_\tQ.\cA.j_\tQ$, we see that $x\in\ker(\br_\tQ)\cap j_\tQ.\cA.j_\tQ$, which is contained in the Jacobson radical $J(j_\tQ.\cA.j_\tQ)$.  Set $v:=j_\tQ+\sum\limits_{n=1}^\infty x^n$,   well-defined by the nilpotency of $x$, and $t_{\varphi^{-1}}:=\widetilde t_{\varphi^{-1}}.v$.  We have
\[
s_\varphi.t_{\varphi^{-1}}=s_\varphi.\widetilde t_{\varphi^{-1}}.\left(j_\tQ+\sum_{n=1}^\infty x^n\right)=(j_\tQ-x).(j_\tQ-x)^{-1}=j_\tQ,
\]
as $j_\tQ$ is the identity of $j_\tQ.\cA.j_\tQ$.  We then compute 
\[
(t_{\varphi^{-1}}.s_\varphi)^2=t_{\varphi^{-1}}.j_\tQ.s_\varphi=t_{\varphi^{-1}}.s_\varphi,
\]
so $t_{\varphi^{-1}}.s_\varphi$ is an idempotent of $i_\tP.\cA.i_\tP$.  It is also nonzero, as
\[
\br_\tP(t_{\varphi^{-1}}.s_\varphi)=\br_\tP(\widetilde t_{\varphi^{-1}}.v.s_\varphi)=\overline i_\tP.u_\varphi^\dag.(\overline j_\tQ)^3.u_\varphi.\overline i_\tP=(\overline i_\tP)^5=\overline i_\tP\neq 0\in\cA(\tP)
\]
as $i_\tP$ was assumed to be local.  The idempotent $i_\tP$ was also assumed to be primitive, so the algebra $i_\tP.\cA.i_\tP$ has a unique nonzero idempotent, namely $i_\tP$ itself.  It follows that that $t_{\varphi^{-1}}.s_\varphi=i_\tP$, as desired.
\end{proof}
\end{prop}

The alternative characterization of the action of $\theta_\varphi$ may seem reminiscent of the discussion of local unital Puig categories in \S\ref{sec:Block_fusion_and_local_categories}.  This is no accident:

\begin{cor}\label{cor:twisted_unit_transport_vs_unital_Puig_iso}
Let $\cA$ be a divisible $\tS$-algebra.  Given $\varphi\in\Iso_{\fF_\tS(\cA)}(\tP,\tQ)$, $\gamma\in\cL\cP(\cA^\tP)$, and $\delta\in\cL\cP(\cA^\tQ)$, we have $\theta_\varphi(\gamma)=\delta$ if and only if $\varphi\in\Iso_{\fL_\tS^\times(\cA)}(\tP_\gamma,\tQ_\delta)$.
\begin{proof}
Given $i_\tP\in\gamma$ and $j_\tQ\in\delta$, each condition is equivalent to the existence of
\[
s_\varphi\in j_\tQ.{^\varphi\!\!\cA^\tP}.i_\tP\qquad\textrm{and}\qquad
t_{\varphi^{-1}}\in i_\tP.{^{\varphi^{-1}}\!\!\cA^\tQ}.j_\tQ
\]
satisfying $i_\tP=t_{\varphi^{-1}}.s_\varphi$ and $j_\tQ=s_\varphi.t_{\varphi^{-1}}$ by Proposition \ref{prop:PRS_algebra_transpotent_lifting} and Lemma \ref{lem:Puig_morphisms_via_transpotents}, respectively.  The result follows.
\end{proof}
\end{cor}

So long as $\cA$ has all twisted units, the isomorphism $\varphi\in\Iso_{\fF_\tS(\cA)}(\tP,\tQ)$ induces not only a bijection between the local points of $\tP$ and $\tQ$, but between the entire pointed Brown posets $s(\tP,\cA)$ and $s(\tQ,\cA)$.   For $\tR\leq \tP$, let $\varphi|_{\tR}\in\Iso_{\fF_\tS(\cA)}(\tR,\varphi \tR)$ be the restriction of $\varphi$.  We then have the bijection
\[
\theta_{\varphi|_{\tR}}:\cL\cP(\cA^{\tR})\to\cL\cP(\cA^{\varphi \tR}),
\]
which exists by the assumption that $\cA(\varphi|_{\tR})^\times\neq\emptyset$.  Piecing together the various $\theta_{\varphi|_\tR}$ as $\tR$ ranges over all subgroups of $\tP$, we obtain
\[
\Theta_\varphi:s(\tP,\cA)\to s(\tQ,\cA):\tR_\varepsilon\mapsto(\varphi \tR)_{\theta_{\varphi|_\tR}(\varepsilon)}.
\]
It is clear that $\Theta_\varphi$ is bijection of sets.  In fact, much more is true:  $\Theta_\varphi$ is an equivariant, multiplicity-preserving poset isomorphism.  Recall that the \emph{multiplicity} of $\tR_\varepsilon\in s(\tP,\cA)$ is the number of elements $m(\tR_\varepsilon)$ in $\varepsilon$ that appear in a primitive idempotent decomposition of $1_{\cA^\tR}$, and for $\tR'_{\varepsilon'}\leq\tR_\varepsilon$ the \emph{relative multiplicity} of $\tR'_{\varepsilon'}$ in $\tR_{\varepsilon}$ is the number of elements $m(\tR'_{\varepsilon'},\tR_\varepsilon)$ in $\varepsilon'$ that appear in a primitive decomposition of any $e\in\varepsilon$, viewed as an idempotent in $\cA^{\tR'}$.

\begin{prop}\label{prop:pointed_Brown_poset_transfer_of_structure}
Let $\cA$ be a divisible $\tS$-algebra that has all twisted units.  For any isomorphism $\varphi\in\Iso_{\fF_\tS(\cA)}(\tP,\tQ)$, local pointed subgroups $\tR'_{\varepsilon'}\leq\tR_\varepsilon\in s(\tP,\cA)$, and group element $\tp\in\tP$, the bijection $\Theta_\varphi:s(\tP,\cA)\to s(\tQ,\cA)$ satisfies:
\begin{enumerate}
\item\label{prop:pointed_Brown_poset_transfer_of_structure_i} $\Theta_\varphi(\tR'_{\varepsilon'})\leq\Theta_\varphi(\tR_\varepsilon)$,
\item\label{prop:pointed_Brown_poset_transfer_of_structure_ii} $\Theta_\varphi({^\tp(\tR_\varepsilon)})={^{\varphi (\tp)}\Theta_\varphi(\tR_\varepsilon)}$,
\item\label{prop:pointed_Brown_poset_transfer_of_structure_iii}  $m(\tR_\varepsilon)=m(\Theta_\varphi(\tR_\varepsilon))$, and
\item\label{prop:pointed_Brown_poset_transfer_of_structure_iv}$m(\tR'_{\varepsilon'},\tR_\varepsilon)=m(\Theta_\varphi(\tR'_{\varepsilon'}),\Theta_\varphi(\tR_\varepsilon))$.
\end{enumerate}
\begin{proof}
Pick $i_\tR\in\varepsilon$, $j_{\varphi \tR}\in\theta_{\varphi|_\tR}(\varepsilon)$, and  elements 
\[
s_\varphi\in j_{\varphi\tR}.{^\varphi\!\!\cA^\tR}.i_\tR\qquad\textrm{and}
\qquad t_{\varphi^{-1}}\in i_\tR.{^{\varphi^{-1}}\!\!\!\cA^{\varphi \tR}}.j_{\varphi\tR}
\]
such that 
$i_\tR=t_{\varphi^{-1}}.s_\varphi$ and $j_{\varphi \tR}=s_\varphi.t_{\varphi^{-1}}$, as provided by Proposition \ref{prop:PRS_algebra_transpotent_lifting}.
\begin{enumerate}
\item We must find an idempotent in $\theta_{\varphi|_{\tR'}}(\varepsilon')$ contained in $j_{\varphi \tR}$.  As $\tR'_{\varepsilon'}\leq \tR_\varepsilon$, we may pick $i'_{\tR'}\in\varepsilon'$ such that $i'_{\tR'}\leq i_\tR$ as idempotents in $\cA^{\tR'}$.  Set
\[
j'_{\varphi \tR'}:=s_\varphi.i'_{\tR'}.t_{\varphi^{-1}}
\]
and compute
\[
\qquad\qquad(j'_{\varphi \tR'})^2=(s_\varphi.i'_{\tR'}.t_{\varphi^{-1}}).(s_\varphi.i'_{\tR'}.t_{\varphi^{-1}})=
s_\varphi.i'_{\tR'}.i_\tR.i'_{\tR'}.t_{\varphi^{-1}}=s_\varphi.i'_{\tR'}.t_{\varphi^{-1}}=j'_{\varphi \tR'},
\]
so that $j'_{\varphi \tR'}$ is an idempotent.  It is clear that $j'_{\varphi \tR'}\in\cA^{\varphi \tR'}$ and the primitivity of $i'_{\tR'}$ implies the same of $j'_{\varphi \tR'}$.  Now set
\[
s'_\varphi:=j'_{\varphi \tR'}.s_\varphi.i'_{\tR'}\qquad\textrm{and}\qquad
t'_{\varphi^{-1}}:= i'_{\tR'}.t_{\varphi^{-1}}.j'_{\varphi \tR'}.
\]
We have
\begin{multline*}
\qquad\qquad t'_{\varphi^{-1}}.s'_\varphi=( i'_{\tR'}.t_{\varphi^{-1}}.j'_{\varphi \tR'}).
(j'_{\varphi \tR'}.s_\varphi.i'_{\tR'})=
i'_{\tR'}.t_{\varphi^{-1}}.j'_{\varphi\tR'}.s_\varphi.i'_{\tR'}
\\
=i'_{\tR'}.t_{\varphi^{-1}}.(s_\varphi.i'_{\tR'}.t_{\varphi^{-1}}).s_\varphi.i'_{\tR'}=i'_{\tR'}.i_\tR.i'_{\tR'}.i_\tR.i'_{\tR'}=i'_{\tR'}
\end{multline*}
and similarly $s'_\varphi.t'_{\varphi^{-1}}=j'_{\varphi \tR'}$.  As $s_\varphi'\in j'_{\varphi\tR'}.{^\varphi\!\!\cA^{\tR'}}.i'_{\tR'}$ and $t'_{\varphi^{-1}}\in i'_{\tR'}.{^{\varphi^{-1}}\!\!\!\cA^{\varphi \tR'}}.j'_{\varepsilon\tR'}$, Proposition \ref{prop:PRS_algebra_transpotent_lifting} gives us $[j'_{\varphi \tR'}]=\theta_{\varphi|_{\tR'}}(\varepsilon')$.  On the other hand, we compute
\[
\qquad\qquad j_{\varphi \tR}.j'_{\varphi \tR'}.j_{\varphi\tR}=
s_\varphi.t_{\varphi^{-1}}.s_\varphi.i'_{\tR'}.t_{\varphi^{-1}}.s_\varphi.t_{\varphi^{-1}}=
s_\varphi.i_\tR.i'_{\tR'}.i_\tR.t_{\varphi^{-1}}=
s_\varphi.i'_{\tR'}.t_{\varphi^{-1}}=j'_{\varphi \tR'},
\]
 so that $j'_{\varphi \tR'}\leq j_{\varphi \tR}$.  The result is proved.

\item We must show that $\theta_{\varphi|_{^\tp\tR}}({^\tp\varepsilon}) ={^{\varphi(\tp)}\theta_{\varphi|_\tR}(\varepsilon)}$.  Observe that $\tp^{-1}\cdot 1_\cA\in\cA^\times\cap{^{c_\tp^{-1}}\!\!\!\cA^{^\tp\tR}}$, and as such determines the twisted  unit $\overline{\tp^{-1}}\in\cA(c_\tp^{-1}|_{^\tp\tR})^\times$.  Similarly, $\varphi(\tp)$ determines the twisted unit $\overline{\varphi(\tp)}\in\cA(c_{\varphi(\tp)}|_{\varphi\tR})^\times$.  Given $u_\varphi\in\cA(\varphi|_\tR)^\times$, it follows from Lemma \ref{lem:twisted_units_properties_i} that we have
\[
\overline{\varphi(\tp)}.u_\varphi.\overline{\tp^{-1}}\in\cA(c_{\varphi(\tp)}|_{\varphi\tR}\circ\varphi|_\tR\circ c_\tp^{-1}|_{^\tp\tR})^\times=\cA(\varphi|_{^\tp\tR})^\times,
\]
as $\varphi$ is defined on all of $\tP$.  Then for $e\in\varepsilon$, we have
\[
\qquad\qquad\theta_{\varphi|_{^\tp\tR}}({^\tp\varepsilon})=\left[\overline{\varphi(\tp)}.u_{\varphi|_\tR}.\overline{\tp^{-1}}.{^\tp e}.\overline \tp.u_{\varphi|_\tR}^\dag.\overline{\varphi(\tp)^{-1}}\right]
={^{\varphi(\tp)}\!\!\left[
u_{\varphi|_\tR}.e.u_{\varphi|_\tR}^\dag\right]}={^{\varphi(\tp)}\theta_{\varphi|_\tR}(\varepsilon)},
\]
as desired.
\item  
Clearly $m(\tR_\varepsilon)$ is equal to the number of elements in $\br_\tR(\varepsilon)$ that appear in a primitive decomposition of $1_{\cA(\tR)}$, and similarly $m(\Theta_\varepsilon(\tR_\varepsilon))$ is the number of elements of  $\br_{\varphi \tR}(\theta_{\varphi|\tR}(\varepsilon))=\overline\theta_{\varphi|_\tR}(\br_\tR(\varepsilon))$ in a primitive decomposition of $1_{\cA(\varphi \tR)}$.  By Lemma \ref{lem:twisted_units_properties_v}, conjugation by any twisted unit in $\cA(\varphi|_\tR)^\times$ induces an algebra isomorphism $\cA(\tR)\to\cA(\varphi \tR)$, which takes $\varepsilon$ to $\overline\theta_{\varphi|_\tR}(\varepsilon)$ by the definition of $\overline\theta_\varphi$.  The result follows.

\item The argument is just an extended version of the previous one.  Let
\[
i_\tR=\sum_{i'_{\tR'}\in I'} i'_{\tR'}
\]
be a primitive idempotent decomposition of $i_\tR$ in $\cA^{\tR'}$, and let $I_0'\subseteq I'$ be the set of idempotents belonging to $\varepsilon'$.  For each $i'_{\tR'}\in I_0'$, the element $s_\varphi.i'_{\tR'}.t_{\varphi^{-1}}$ is an idempotent in $\theta_{\varphi|_{\tR'}}(\varepsilon')$ contained in some idempotent of $\theta_{\varphi|_\tR}(\varepsilon)$, as seen in  \ref{prop:pointed_Brown_poset_transfer_of_structure_i}.  Moreover, if $k'_{\tR'}$ is a different element of $I_0'$, we have
\[
s_\varphi.i'_{\tR'}.t_{\varphi^{-1}}.s_\varphi.k'_{\tR'}.t_{\varphi^{-1}}=s_\varphi.i'_{\tR'}.i_R.k'_{\tR'}.t_{\varphi^{-1}}=s_\varphi.i'_{\tR'}.k'_{\tR'}.t_{\varphi^{-1}}=0,
\]
so $\{s_\varphi.i'_{\tR'}.t_{\varphi^{-1}}\ \big|\ i'_{\tR'}\in I_0'\}$ is a set of orthogonal idempotents in $\theta_{\varphi|_{\tR'}}(\varepsilon')$, each contained in an idempotent of $\theta_{\varphi|_\tR}(\varepsilon)$.  Thus $|I_0'|=m(\tR'_{\varepsilon'},\tR_{\varepsilon})\leq m(\Theta_\varphi(\tR'_{\varepsilon'}),\Theta_\varphi(\tR_\varepsilon))$.  The situation is symmetric, so running the same argument with $\varphi^{-1}$ in place of $\varphi$ yields the opposite inequality, and the result is proved.
\end{enumerate}
\end{proof}
\end{prop}

Our last ingredient is the notion of \emph{local invariant decomposition}, which we briefly recall.  Returning to the full action of $\tS$ on $\cA$, we consider an idempotent decomposition
\[
1_\cA=\sum_{i\in I}i
\]
such that $^\ts i=\ts\cdot i\cdot \ts^{-1}\in I$ for all $\ts\in \tS$ and $i\in I$.  Such a decomposition is called \emph{invariant}, and $I$ itself may be viewed as an $\tS$-set in this case.  In particular, each idempotent $i$ has an $\tS$-stablizer $\tS_i$, and if $i$ is a primitive local idempotent in $\cA^{\tS_i}$ for all $i\in I$, we say the invariant decomposition is \emph{local}.

The key facts we will need were proved in \cite{PuigSurUneTheoremeDeGreen} and are summarized in \cite[\S24]{ThevenazBook}:
\begin{enumerate}
\item Local invariant decompositions always exist (so long as $\tS$ is a $p$-group), and are unique up to $(\cA^\tS)^\times$-conjugacy.
\item If $I$ is a local invariant decomposition and $i^+:=\tra_{\tS_i}^\tS i$ is the sum of the elements in the $\tS$-orbit of $i\in I$, then $i^+$ is a primitive idempotent of $\cA^\tS$.  In particular, the collection of all such orbit-sums is a primitive idempotent decomposition of $\cA^\tS$.
\end{enumerate}

We can now prove the main result of this section.

\begin{theorem}\label{thm:main_ii_implies_i}
Let $\cA$ be a divisble $\tS$-algebra that has all twisted units.  Then for all $\varphi\in\Iso_{\fF_\tS(\cA)}(\tP,\tQ)$, we have $\cA^\times\cap{^\varphi\!\!\cA^\tP}\neq\emptyset$.
\begin{proof}
Fix local invariant decompositions
\[
1_\cA=\sum_{e\in E} e\qquad\textrm{and}\qquad 1_\cA=\sum_{f\in F} f
\]
of $\cA$ viewed as a $\tP$- and $\tQ$-algebra, respectively.  These decompositions may be partitioned according to the conjugacy class of each idempotent's stabilizer and the local point of that stabilizer determined by the idempotent:  For $\tR_\varepsilon\in s(\tP,\cA)$ and $\tT_\zeta\in s(\tQ,\cA)$, set
\[
E_{\tR_\varepsilon}:=\left\{ e\in E\ \big|\ (\tP_e)_{[e]}\cong_\tP \tR_\varepsilon\right\}\qquad\textrm{and}\qquad
F_{\tT_\zeta}:=\left\{ f\in F\ \big|\ (\tQ_f)_{[f]}\cong_\tQ \tT_\zeta\right\}.
\]
(The authors apologize for the awkward notation.)  Our first goal is to establish a bijection $\sigma_{\tR_\varepsilon}:E_{\tR_\varepsilon}\to F_{\Theta_\varphi(\tR_\varepsilon)}$ for each $\tR_{\varepsilon}\in s(\tP,\cA)$.  We will do this by computing the multiplicities $m(\tR_\varepsilon)$ and $m(\Theta_\varphi(\tR_\varepsilon))$, which are known to be equal by Proposition \ref{prop:pointed_Brown_poset_transfer_of_structure}\ref{prop:pointed_Brown_poset_transfer_of_structure_iii}, in terms of $E$ and $F$.  Clearly we have
\[
m(\tR_\varepsilon)=\sum_{e\in E}m(\tR_\varepsilon,(\tP_e)_{[e]})\qquad\textrm{and}\qquad m(\Theta_\varphi(\tR_\varepsilon))=\sum_{f\in F}m(\Theta_\varphi(\tR_\varepsilon),(\tQ_f)_{[f]}),
\]
but, as $m(\tR_\varepsilon,(\tP_e)_{[e]})$ need not equal $m(\tR_\varepsilon,{^\tp((\tP_e)_{[e]})})$ for $\tp\in \tP$, we cannot directly reexpress $m(\tR_\varepsilon)$ in terms of our partion $E_{\tR_\varepsilon}$.  However, we can get around this obstacle by considering the idempotent defined by adding together the $\tP$-orbit of $e$, and considering the multiplicity of $\tR_\varepsilon$ there:
\[
m(\tR_\varepsilon,\tP_{[e^+]}):=\sum_{\tp\in [\tP/\tP_e]} m(\tR_\varepsilon,{^\tp((\tP_e)_{[e]})}).
\]
Observe that $m(\tR_\varepsilon,\tP_{[e^+]})=m(\tR_\varepsilon,\tP_{[\tilde{e}^+]})$ so long as $e$ and $\tilde{e}$ both live in the same $E_{\widehat{\tR}_{\widehat{\varepsilon}}}$, even if they do not live in the same $\tP$-orbit, and both relative multiplicities are nonzero if and only if $\tR_\varepsilon\leq_\tP\widehat{\tR}_{\widehat{\varepsilon}}$.  Write $m^+(\tR_\varepsilon,\widehat \tR_{\widehat\varepsilon})$ for this common value. Thus, letting $[\tP\backslash s(\tP,\cA)]$ be a set of representatives of the $\tP$-conjugacy classes of the pointed Brown poset, and noting that there are $[\tP:\tP_e]$ elements of $E_{(\tP_e)_{[e]}}$ that contribute to $e^+$, we have
\[
m(\tR_\varepsilon)=
\sum_{
\substack{
{\widehat{\tR}}_{\widehat{\varepsilon}}\in [\tP\backslash s(\tP,\cA)] \\
\tR_\varepsilon\leq_\tP\widehat{\tR}_{\widehat{\varepsilon}}
}}
\frac{\big|E_{\widehat{\tR}_{\widehat{\varepsilon}}}\big|}{[\tP:\widehat{\tR}]}\cdot m^+(\tR_\varepsilon,\widehat{\tR}_{\widehat{\varepsilon}})
\]
By isolating the term indexed by $\tR_\varepsilon$  and rearranging, we have
\[
\big|E_{\tR_\varepsilon}\big|=\frac{[\tP:\tR]}{m^+(\tR_\varepsilon,\tR_\varepsilon)}\cdot\left(m(\tR_\varepsilon)-
\sum_{
\substack{
{\widehat{\tR}}_{\widehat{\varepsilon}}\in [\tP\backslash s(\tP,\cA)] \\
\tR_\varepsilon\lneq_\tP\widehat{\tR}_{\widehat{\varepsilon}}
}}
\frac{\big|E_{\widehat{\tR}_{\widehat{\varepsilon}}}\big|}{[\tP:\widehat{\tR}]}\cdot m^+(\tR_\varepsilon,\widehat{\tR}_{\widehat{\varepsilon}})
\right).
\]
The same considerations in $F$ yield 
\[
\big|F_{\Theta_\varphi(\tR_\varepsilon)}\big|=\frac{[\tQ:\varphi\tR]}{m^+(\Theta_\varphi(\tR_\varepsilon),\Theta_\varphi(\tR_\varepsilon))}\cdot\left(m(\Theta_\varphi(\tR_\varepsilon))-
\sum_{
\substack{
{\widehat{\tT}}_{\widehat{\zeta}}\in [\tQ\backslash s(\tQ,\cA)] \\
\Theta_\varphi(\tR_\varepsilon)\lneq_\tQ\widehat{\tT}_{\widehat{\zeta}}
}}
\frac{\big|F_{\widehat{\tT}_{\widehat{\zeta}}}\big|}{[\tQ:\widehat{\tT}]}\cdot m^+(\Theta_\varphi(\tR_\varepsilon),\widehat{\tT}_{\widehat{\zeta}})
\right).
\]
This suggests the equality $|E_{\tR_\varepsilon}|=|F_{\Theta_\varphi(\tR_\varepsilon)}|$ may be proved by a downward induction on $[\tP\backslash s(\tP,\cA)]$ (ordered by $\tP$-subconjugacy) if all the corresponding terms in these expressions can be equated.  All of the relevant indices of subgroups are obviously equal ($\tQ=\varphi\tP$ so $[\tQ:\varphi\tR]=[\tP: \tR]$, etc.) as $\varphi$ is an injection.  That $m(\Theta_\varphi(\tR_\varepsilon))=m(\tR_\varepsilon)$ is Proposition \ref{prop:pointed_Brown_poset_transfer_of_structure}\ref{prop:pointed_Brown_poset_transfer_of_structure_iii}.  That $|F_{\Theta_\varphi(\widehat{\tR}_{\widehat{\varepsilon}})}|=|E_{\widehat{\tR}_{\widehat{\varepsilon}}}|$ for all $\widehat{\tR}_{\widehat{\varepsilon}}\gneq \tR_\varepsilon$ in $[\tP\backslash s(\tP,\cA)]$ is the inductive hypothesis, whose base case (when $\tR=\tP$) follows from the observation that $|E_{\tP_\gamma}|=m(\tP_\gamma)$ and $|F_{\Theta_\varphi(\tP_\gamma)}|=m(\Theta_\varphi(\tP_\gamma)$, and another appliction of Proposition \ref{prop:pointed_Brown_poset_transfer_of_structure}\ref{prop:pointed_Brown_poset_transfer_of_structure_iii}.

We are left to prove that $m^+(\tR_\varepsilon,\widehat{\tR}_{\widehat{\varepsilon}})=m^+(\Theta_\varphi(\tR_\varepsilon),\Theta_\varphi(\widehat{\tR}_{\widehat{\varepsilon}}))$ for all $\widehat{\tR}_{\widehat{\varepsilon}}\geq \tR_\varepsilon$.  This is a consequence of the equivariance and relative multiplicity preservation of $\Theta_\varphi$, as
\begin{multline*}
m^+(\tR_\varepsilon,\widehat{\tR}_{\widehat{\varepsilon}})=\sum_{\tp\in [\tP/\widehat{\tR}]}m(\tR_\varepsilon,{^\tp(\widehat{\tR}_{\widehat{\varepsilon}})})
=\sum_{\tp\in [\tP/\widehat{\tR}]}m(\Theta_\varphi(\tR_\varepsilon),\Theta_\varepsilon({^\tp(\widehat{\tR}_{\widehat{\varepsilon}}})))\\
=
\sum_{\varphi(\tp)\in [\tP/\varphi(\widehat{\tR})]}m(\Theta_\varphi(\tR_\varepsilon),{^{\varphi(\tp)}\Theta_\varphi(\widehat{\tR}_{\widehat{\varepsilon}}}))=
m^+(\Theta_\varphi(\tR_\varepsilon),\Theta_\varphi(\widehat{\tR}_{\widehat{\varepsilon}}))
\end{multline*}
by Proposition \ref{prop:pointed_Brown_poset_transfer_of_structure}\ref{prop:pointed_Brown_poset_transfer_of_structure_iv} and \ref{prop:pointed_Brown_poset_transfer_of_structure_ii}.  This completes the necessary data for our inductive argument, and we have $|E_{\tR_\varepsilon}|=|F_{\Theta_\varphi(\tR_\varepsilon)}|$ for all $\tR_\varepsilon\in s(\tP,\cA)$.

At this point we have established the existence of bijections $\sigma_{\tR_\varepsilon}:E_{\tR_\varepsilon}\to F_{\Theta_\varepsilon(\tR_\varepsilon)}$ for all $\tR_\varepsilon\in [\tP\backslash s(\tP,\cA)]$.  As $E_{\tR_\varepsilon}$ is, as a $\tP$-set, a disjoint union of copies of $[\tP/\tR]$, and $E_{\Theta_\varepsilon(\tR_\varepsilon)}$ is likewise a disjoint union of copies of $[\tQ/\varphi \tR]$, we ay choose $\sigma_{\tR_{\varepsilon}}$ to be equivariant, i.e., $\sigma_{\tR_\varepsilon}({^\tp e})={^{\varphi(\tp)}\sigma_{\tR_\varepsilon}(e)}$ for all $\tp\in \tP$ and $e\in E_{\tR_\varepsilon}$.  We  further assume that $\theta_\varepsilon([e])=[\sigma_{\tR_\varepsilon}(e)]$ for all $e\in E_{\tR_\varepsilon}$.

For each $e\in E_{\tR_\varepsilon}$, pick elements $s_\varphi^e\in\sigma_{\tR_\varepsilon}(e).{^\varphi\!\!\cA^\tR}.e$ and $t_{\varphi^{-1}}^e\in e.{^{\varphi^{-1}}\!\!\!\cA^{\varphi \tR}}.\sigma_{\tR_\varepsilon}(e)$ such that $e=t_{\varphi^{-1}}.s_\varphi^e$ and $\sigma_{\tR_\varepsilon}(e)=s_\varphi^e.t_{\varphi^{-1}}^e$, as per Proposition \ref{prop:PRS_algebra_transpotent_lifting}.  Setting
\[
u_\varphi^e:=\tra_{\Delta(\varphi|_\tR,\tR)}^{\Delta(\varphi,\tP)}(s_\varphi^e)\qquad\textrm{and}\qquad
v_{\varphi^{-1}}^e:=\tra_{\Delta(\varphi^{-1}|_{\varphi\tR},\varphi \tR)}^{\Delta(\varphi^{-1},\tQ)}(t_{\varphi^{-1}}^e)
\]
yields elements that satisfy $e^+=v_{\varphi^{-1}}^e.u_\varphi^e$ and $\sigma_{\tR_\varepsilon}(e)^+=u_\varepsilon^e.v_{\varphi^{-1}}^e$.  Therefore
\[
u:=\sum_{\substack{\tR_{\varepsilon}\in [\tP\backslash s(P,\cA)]\\ e\in [P\backslash E_{\tR_\varepsilon}]}}u_\varphi^e\qquad\textrm{and}\qquad
v:=\sum_{\substack{\tR_{\varepsilon}\in [\tP\backslash s(\tP,\cA)]\\ e\in [\tP\backslash E_{\tR_\varepsilon}]}}v_{\varphi^{-1}}^e
\]
are mutual inverses in $\cA$.  As $u\in{^\varphi\!\!\cA^\tR}$ by construction, the result is proved.
\end{proof}
\end{theorem}

Summarizing these results, we have

\begin{cor}\label{cor:1_2_equiv}
Let $\cA$ be a divisible $\tS$-algebra.  Then $\cA$ possesses a unital $(\tS,\tS)$-invariant $\cO$-basis if and only if $\cA$ has all twisted units.
\begin{proof}
Then ``only if'' implication was spelled out in the introduction to this section:  If $X$ is a unital $(\tS,\tS)$-invariant $\cO$-basis and $\varphi\in\Iso_{\fF_\tS(\cA)}(P,Q)$, then $\varphi$ is realized by some $x\in{^\varphi\! X^\tP}$, which, being a unit, yields the twisted unit $\br_\varphi(x)\in\cA(\varphi)^\times$.

On the other hand, if $\cA$ has all twisted units, then for any $\varphi\in\Iso_ {\fF_\tS(\cA)}$ we have $\cA^\times\cap{^\varphi\!\!\cA^\tP}\neq\emptyset$ by Theorem \ref{thm:main_ii_implies_i}.  Proposition \ref{prop:unit_stabilizers_implies_unit_basis} then finishes the proof.
\end{proof}
\end{cor}

When we take $\cA=\cS$, this proves the equivalence of Conjectures $\ref{conj:main_strong_i}$ and $\ref{conj:main_strong_ii}$.

\section{Balanced corner algebras}\label{sec:balanced_algebras}

Let $\cA$ be a divisible $\tS$-algebra.  In this section we introduce the final structural property equivalent to $\cA$'s having a unital $(\tS,\tS)$-invariant $\cO$-basis. Our method will grow out of the work of the previous section, so we summarize the main ingredients that went into the proof of Theorem \ref{thm:main_ii_implies_i}:
\begin{enumerate}
\item Any $\varphi\in\Iso_{\fF_\tS(\cA)}(\tP,\tQ)$ determines a bijection $\theta_\varphi:\cL\cP(\cA^\tP)\to\cL\cP(\cA^\tQ)$.
\item This bijection preserves the multiplicities of local points:  $m(\tP_\gamma)=m(\tQ_{\theta_\varphi(\gamma)})$.
\item For  $i_\tP\in\gamma\in\cL\cP(\cA^\tP)$ and $j_\tQ\in\theta_\varphi(\gamma)\in\cL\cP(\cA^\tQ)$  there are transporting elements $s_\varphi\in j_\tQ.{^\varphi\!\!\cA^\tP}.i_\tP$ and $t_{\varphi^{-1}}\in i_\tP.{^{\varphi^{-1}}\!\!\!\cA^\tQ}.j_\tQ$ such that $i_\tP=t_{\varphi^{-1}}.s_\varphi$ and $j_\tQ=s_\varphi.t_{\varphi^{-1}}$.
\end{enumerate}
While before these facts were derived from properties of twisted units, now we take them as the starting point for our discussion.

Unlike our earlier properties, which dealt solely with the structure of $\cA$, we will now work in the context of an ambient interior $\tS$-algebra $\widehat{\cA}$, in which $\cA$ is embedded as a corner algebra $\cA=\ell.\widehat{\cA}.\ell$ for some idempotent $\ell\in\widehat{\cA}^\tS$.  We assume that $\widehat\cA$ possesses a unital $(\tS,\tS)$-invariant $\cO$-basis, and seek a condition on the choice of $\ell$ that will guarantee that this property is inherited by $\cA$. 

The assumption on $\widehat{\cA}$ implies that it has all twisted units by Corollary \ref{cor:1_2_equiv}, so for any $\varphi\in\Iso_{\fF_\tS(\widehat{\cA})}(\tP,\tQ)$ we have the induced bijection of local points
\[
\theta_\varphi:\cL\cP(\widehat \cA^\tP)\to\cL\cP(\widehat \cA^\tQ).
\]

\begin{definition}
Let $\widehat{\cA}$ be a divisible  $\tS$-algebra with unital $(\tS,\tS)$-invariant $\cO$-basis, $\ell$ an idempotent in $\widehat{\cA}^\tS$, and $\cA:=\ell.\widehat{\cA}.\ell$ the corresponding corner algebra.   Then $\cA$ is \emph{balanced} in $\widehat{\cA}$ if for all $\varphi\in\Iso_{\fF_\tS(\cA)}(\tP,\tQ)$ and all local points $\widehat \gamma\in\cL\cP(\widehat{\cA}^\tP)$ we have 
\[
m(\tP_{\widehat{\gamma}},\tS_{[\ell]})=m(\tQ_{\theta_\varphi(\widehat{\gamma})},\tS_{[\ell]}).
\]
\end{definition}

Note that the balance condition tests only against the isomorphisms of $\fF_\tS(\cA)$, not those of the larger $\fF_\tS(\widehat{\cA})$.  Also note that the condition implies that for all $\tP_{\widehat{\gamma}}$ in $s(\tS,\widehat{\cA})$ such that $\gamma:=\widehat{\gamma}\cap\cA\neq\emptyset$, if there is some
\[
\varphi\in\Iso_{\cL_\tS^\times(\widehat{\cA})}(\tP_{\widehat{\gamma}},\tQ_{\widehat{\delta}})\cap\Iso_{\fF_\tS(\cA)}(\tP,\tQ),
\]
then $\delta:=\widehat{\delta}\cap\cA\neq\emptyset$ as well. 

We can reformulate the balance condition without  reference to an ambient algebra:

\begin{definition}
The divisible $\tS$-algebra $\cA$  is \emph{intrinsically balanced} if
\begin{enumerate}
\item  All $\fF_\tS(\cA)$-isomorphisms are $\cA^\times$-isofusions for all choices of source pointed group:  For every  $\varphi\in\Iso_{\fF_\tS(\cA)}(\tP,\tQ)$ and local point $\gamma\in\cL\cP(\cA^\tP)$, there is a (unique) local point $\delta\in\cL\cP(\cA^{\tQ})$ such that $\varphi\in\Iso_{\fL_\tS^\times(\cA)}(\tP_\gamma,\tQ_\delta)$.
\item Isomorphic local pointed subgroups of $\cA$ have equal multiplicities:  If $\tP_\gamma$ and $\tQ_\delta$ are isomorphic in $\fL_\tS^\times(\cA)$, then  $m(\tP_\gamma)=m(\tQ_\delta)$.
\end{enumerate}
\end{definition}

\begin{prop}\label{prop:balance_def_equiv}
Let $\widehat\cA$ be a divisble   $\tS$-algebra with a unital $(\tS,\tS)$-invariant $\cO$-basis, $\ell\in\widehat\cA^\tS$ an idempotent, and $\cA:=\ell.\widehat{\cA}.\ell$ the corresponding corner algebra.  Then $\cA$ is balanced in $\widehat{\cA}$ if and only if $\cA$ is intrinsically balanced.
\begin{proof}
First suppose that $\cA$ is balanced in $\widehat\cA$.  Fix $\tP_\gamma\in s(\tS,\cA)$ and let $\tP_{\widehat{\gamma}}$ be the corresponding local pointed group on $\widehat{\cA}$.  If $\varphi\in\Iso_{\fF_\tS(\cA)}(\tP,\tQ)$, we clearly have $\varphi\in\Iso_{\fF_\tS(\widehat \cA)}(\tP,\tQ)$, so by Proposition \ref{prop:unital_fusion_Puig_lifting} we have  $\varphi\in\Iso_{\fL_\tS^\times(\widehat{\cA})}(\tP_{\widehat{\gamma}},\tQ_{\widehat{\delta}})$ for a unique $\widehat{\delta}\in\cL\cP(\widehat{\cA}^\tQ)$. Thus $m(\tQ_{\widehat{\delta}},\tS_{[\ell]})=m(\tP_{\widehat{\gamma}},\tS_{[\ell]})\neq 0$. In particular, $\delta:=\widehat{\delta}\cap\cA$ is nonempty, and so $\tQ_\delta$ is a local pointed group on $\cA$.  Proposition \ref{prop:Puig_categories_of_corner_algebras_are_full} then implies $\varphi\in\Iso_{\fL_\tS^\times(\cA)}(\tP_\gamma,\tQ_\delta)$, and the first condition of the intrinsic balance of $\cA$ is met.  The second condition follows immediately from the observation that the relative multiplicity $m(\tP_{\widehat{\gamma}},\tS_{[\ell]})$, taken in the algebra $\widehat{\cA}$, is equal to the absolute multiplicity $m(\tP_\gamma)$ taken in ${\cA}$, and similarly for $m(\tQ_{\widehat{\delta}},\tS_{[\ell]})$ and $m(\tQ_\delta)$.  Thus $\cA$ is intrinsically balanced.

If  $\cA$ is intrinsically balanced, suppose we are given $\varphi\in\Iso_{\fF_\tS(\cA)}(\tP,\tQ)$ and $\widehat{\gamma}\in\cL\cP(\widehat{\cA}^\tP)$.  If $\gamma:=\widehat{\gamma}\cap\cA\neq\emptyset$, then by hypothesis we have $\varphi\in\Iso_{\fL_\tS^\times}(\tP_\gamma,\tQ_\delta)$ for a unique $\delta\in\cL\cP(\cA^\tQ)$, and $m(\tP_\gamma)=m(\tP_\delta)$.  As $\fL_\tS^\times(\cA)$ embeds fully faithfully  in $\fL_\tS^\times(\widehat{\cA})$ by Proposition \ref{prop:Puig_categories_of_corner_algebras_are_full}, we have $\varphi\in\Iso_{\fL_\tS^\times(\widehat{\cA})}(\tP_{\widehat{\gamma}},\tQ_{\widehat{\delta}})$ for $\widehat{\delta}\in\cL\cP(\widehat{\cA}^\tQ)$ the unique local point containing $\delta$.  As above, we can equate $m(\tP_\gamma)=m(\tP_{\widehat{\gamma}},\tS_{[\ell]})$ and $m(\tQ_\delta)= m(\tQ_{\widehat{\delta}},\tS_{[\ell]})$, we must just show that $\widehat{\delta}=\theta_\varphi(\widehat{\gamma})$.  This is precisely the content of Corollary \ref{cor:twisted_unit_transport_vs_unital_Puig_iso}.  

It remains to consider the case $\widehat{\gamma}\cap\cA=\emptyset$, so that $m(\tP_{\widehat{\gamma}},\tS_{[\ell]})=0$.  We must then have $m(\tQ_{\theta_\varphi(\widehat{\gamma})},\tS_{[\ell]})=0$ as well, or else reach a contradiciton by applying the above argument  to $\varphi^{-1}\in\Iso_{\fF_\tS(\cA)}(\tQ,\tP)$ and $\theta_\varphi(\widehat{\gamma})\in\cL\cP(\widehat{\cA}^\tQ)$.  This proves the result.
\end{proof}
\end{prop}

We are now ready to state and prove our final theorem:

\begin{theorem}\label{thm:balance_and_unital_basis}
Let $\widehat{\cA}$ be a divisible  $\tS$-algebra possessing a unital $(\tS,\tS)$-invariant $\cO$-basis, $\ell\in\cA^\tS$ an idempotent, and $\cA:=\ell. \widehat{\cA}.\ell$ the corresponding corner algebra.  Then $\cA$ possesses a unital $(\tS,\tS)$-invariant $\cO$-basis if and only if $\cA$ is balanced in $\widehat{\cA}$.
\begin{proof}
By Proposition \ref{prop:balance_def_equiv}, we may replace the condition that $\cA$ be balanced in $\widehat{\cA}$ with $\cA$'s being intrinsically balanced, which we do without further comment.

First suppose that $\cA$ has a unital $(\tS,\tS)$-invariant $\cO$-basis $X$.  Let $\varphi\in\Iso_{\fF_\tS(\cA)}(\tP,\tQ)$ be realized by $x\in{^\varphi\! X^\tP}$.  For any $\gamma\in\cL\cP(\cA^\tP)$, we have $\delta:={^x\gamma}\in\cL\cP(\cA^\tQ)$.  One sees immediately that the unit $x\in\cA^\times$ realizes $\varphi:\tP_\gamma\xrightarrow\cong\tQ_\delta$ as an $\cA^\times$-isofusion.   Thus the first condition of intrinsic balance is satisfied.

Consider $\tP_\gamma,\tQ_\delta\in s(\tS,\cA)$ and $\varphi:\tP_\gamma\xrightarrow\cong\tQ_\delta$ an $\cA^\times$-isofusion. By Lemma \ref{lem:Puig_fusion_containment} we have $\Iso_{\fL_\tS^\times(\cA)}(\tP,\tQ)\subseteq\Iso_{\fF_\tS(\cA)}(\tP,\tQ)$, so $\varphi$ is realized by some $x\in{^\varphi\! X^\tP}$.  Then $\delta={^x\gamma}$, and if
$
1_{\cA^\tP}=\sum\limits_{i\in I}i
$
is a primitive idempotent decomposition of $\cA^\tP$, we have $1_{\cA^\tQ}=\sum\limits_{i\in I}{^xi}$ is a primitive idempotent decomposition in $\cA^\tQ$.  Moreover, $i\in I\cap\gamma$ if and only if $^xi\in\delta$.  This establishes $m(\tP_\gamma)=m(\tQ_\delta)$, and so  $\cA$ is intrinsically balanced.

Now suppose that $\cA$ is intrinsically balanced.  By Corollary \ref{cor:1_2_equiv}, it is enough to show that $\cA$ has all twisted units.  In fact, the three key ingredients listed at the beginning of this section are equivalent to $\cA$'s having all twisted units, as we now show:

For each $\varphi\in\Iso_{\fF_\tS(\cA)}(\tP,\tQ)$ let $\theta_\varphi:\cL\cP(\cA^\tP)\to\cL\cP(\cA^\tQ)$ be a bijection of local points that preserves multiplicities and for which there exist transporting elements as  above.  Choose primitive idempotent decompositions
\[
1_{\cA^\tP}=\sum_{i\in I}i\qquad\textrm{and}\qquad 1_{\cA^\tQ}=\sum_{j\in J}j
\]
of $\cA^\tP$ and $\cA^\tQ$, respectively.  If we take $\cL(I)\subseteq I$  to be the subset of local primitive idempotents, and similarly  $\cL(J)\subseteq J$,  we have
\[
1_{\cA(\tP)}=\sum_{i\in\cL (I)}\br_\tP(i)\qquad\textrm{and}\qquad 1_{\cA(\tQ)}=\sum_{j\in\cL (J)}\br_\tQ(j)
\]
as primitive idempotent decompositions of $\cA(\tP)$ and $\cA(\tQ)$.  The assumption that the bijection ${\theta_\varphi:\cL\cP(\cA^\tP)\to\cL\cP(\cA^\tQ)}$ preserves multiplicities gives us a bijection $\vartheta:\cL(I)\to\cL(J)$ such that $[\vartheta(i)]=\theta_\varphi([i])$.  The assumption on the existence of transporting elements allows us to choose, for each $i\in\cL(I)$, elements
\[
s_\varphi^i\in\vartheta(i).{^\varphi\!\!\cA^\tP}.i\qquad\textrm{and}\qquad
t_{\varphi^{-1}}^i\in i.{^{\varphi^{-1}}\!\!\!\cA^\tQ}.\vartheta(i)
\]
such that $i=t_{\varphi^{-1}}^i.s_\varphi^i$ and $\vartheta(i)=s_\varphi^i.t_{\varphi^{-1}}^i$.  We then have $\sum\limits_{i\in\cL(I)}\br_\varphi( s_\varphi^i)$ is a twisted unit in $\cA(\varphi)$ with twisted inverse $\sum\limits_{i\in\cL(I)}\br_{\varphi^{-1}}(t_{\varphi^{-1}}^i)$, and the initial claim is proved.

It therefore suffices to prove that $\cA$'s being intrinsically balanced implies the existence of the bijections $\theta_\varphi$ as above.  This is essentially a restatement of the intrinsic balance condition:  For $\varphi\in\Iso_{\fF_\tS(\cA)}(\tP,\tQ)$  and $\gamma\in\cL\cP(\cA^\tP)$, let $\theta_\varphi(\gamma)\in\cL\cP(\cA^\tQ)$ denote the unique (by Lemma \ref{lem:Puig_category_unique_target}) local point such that $\varphi\in\Iso_{\fL_\tS^\times(\cA)}(\tP_\gamma,\tQ_\delta)$.  As $\cA$ is divisible, this clearly defines a bijection $\cL\cP(\cA^\tP)\to\cL\cP(\cA^\tQ)$ with inverse $\theta_{\varphi}^{-1}=\theta_{\varphi^{-1}}$ defined in the obvious manner.  That $\theta_\varphi$ preserves multiplicities is precisely the second point in the definition of intrinsic balance.  Finally, the existence of transporting elements follows from Lemma \ref{lem:Puig_morphisms_via_transpotents}.  This completes the proof that $\cA$ has all twisted units, and hence an $(\tS,\tS)$-invariant unital $\cO$-basis. 
\end{proof}
\end{theorem}

 Theorem \ref{thm:balance_and_unital_basis} applied to the source algebra $\cS$  yields the last equivalence of Theorem \ref{thm:main}.


\begin{thebibliography}{AKO}

\bibitem[Alp]{AlperinFusion}
J.~L. Alperin.
\newblock Sylow intersections and fusion.
\newblock {\em J. Algebra}, 6:222--241, 1967.

\bibitem[AB]{AlperinBroueBlockFusion}
J.~Alperin and Michel Brou{\'e}.
\newblock Local methods in block theory.
\newblock {\em Ann. of Math. (2)}, 110(1):143--157, 1979.

\bibitem[AKO]{AKO}
Michael Aschbacher, Radha Kessar, and Bob Oliver.
\newblock {\em Fusion systems in algebra and topology}, volume 391 of {\em
  London Mathematical Society Lecture Note Series}.
\newblock Cambridge University Press, Cambridge, 2011.

\bibitem[BLO]{BLO2}
Carles Broto, Ran Levi, and Bob Oliver.
\newblock The homotopy theory of fusion systems.
\newblock {\em J. Amer. Math. Soc.}, 16(4):779--856 (electronic), 2003.

\bibitem[BP]{BrouePuigLocalStructure}
Michel Brou\'{e} and Llu\'{i}s Puig.
\newblock Characters and local structure in {$G$}-algebras.
\newblock {\em J. Algebra}, 63(2):306--317, 1980.

\bibitem[Bra]{BrauerNotesOnRepresentations}
Richard Brauer.
\newblock Notes on representations of finite groups. {I}.
\newblock {\em J. London Math. Soc. (2)}, 13(1):162--166, 1976.

\bibitem[Gel]{GelvinCharacteristicBlockBasisOnline}
Matthew J.~K. Gelvin.
\newblock An observation on the module structure of block algebras.
\newblock {\em Communications in Algebra}, 0(0):1--8, 2019.
\newblock doi: 10.1080/00927872.2019.1617874.

\bibitem[GR]{GelvinReehMinimalCharacteristicBisets}
Matthew Gelvin and Sune~Precht Reeh.
\newblock Minimal characteristic bisets for fusion systems.
\newblock {\em J. Algebra}, 427:345--374, 2015.

\bibitem[GRY]{GelvinReehYalcin}
Matthew Gelvin, Sune~Precht Reeh, and Erg{\"u}n Yal{\c{c}}{\i}n.
\newblock On the basis of the {B}urnside ring of a fusion system.
\newblock {\em J. Algebra}, 423:767--797, 2015.

\bibitem[Gre]{GreenBlocksOfModularRepresentations}
James~A. Green.
\newblock Blocks of modular representations.
\newblock {\em Math. Z.}, 79:100--115, 1962.

\bibitem[Lin1]{LinckelmannOnSpendidDerived}
Markus Linckelmann.
\newblock On splendid derived and stable equivalences between blocks of finite
  groups.
\newblock {\em J. Algebra}, 242(2):819--843, 2001.

\bibitem[Lin2]{LinckelmannBlocksOfMinimalDimension}
Markus Linckelmann.
\newblock Blocks of minimal dimension.
\newblock {\em Arch. Math. (Basel)}, 89(4):311--314, 2007.

\bibitem[Lin3]{LinckelmannTrivialSource}
Markus Linckelmann.
\newblock Trivial source bimodule rings for blocks and {$p$}-permutation
  equivalences.
\newblock {\em Trans. Amer. Math. Soc.}, 361(3):1279--1316, 2009.

\bibitem[Lin4]{LinckelmannBookI}
Markus Linckelmann.
\newblock {\em The block theory of finite group algebras. {V}ol. {I}},
  volume~91 of {\em London Mathematical Society Student Texts}.
\newblock Cambridge University Press, Cambridge, 2018.

\bibitem[Lin5]{LinckelmannBookII}
Markus Linckelmann.
\newblock {\em The block theory of finite group algebras. {V}ol. {II}},
  volume~92 of {\em London Mathematical Society Student Texts}.
\newblock Cambridge University Press, Cambridge, 2018.

\bibitem[Pui1]{PuigSurUneTheoremeDeGreen}
Llu\'{i}s Puig.
\newblock Sur un th\'{e}or\`eme de {G}reen.
\newblock {\em Math. Z.}, 166(2):117--129, 1979.

\bibitem[Pui2]{PuigLocalFusions}
Llu{\'i}s Puig.
\newblock Local fusions in block source algebras.
\newblock {\em J. Algebra}, 104(2):358--369, 1986.

\bibitem[Pui3]{PuigConstructionOfModules}
Llu{\'{\i}}s Puig.
\newblock Pointed groups and construction of modules.
\newblock {\em J. Algebra}, 116(1):7--129, 1988.

\bibitem[Pui4]{PuigFrobeniusCategories}
Llu{\' i}s Puig.
\newblock Frobenius categories.
\newblock {\em J. Algebra}, 303(1):309--357, 2006.

\bibitem[Pui5]{PuigBook}
Llu{\'{\i}}s Puig.
\newblock {\em Frobenius categories versus {B}rauer blocks}, volume 274 of {\em
  Progress in Mathematics}.
\newblock Birkh\"auser Verlag, Basel, 2009.
\newblock The Grothendieck group of the Frobenius category of a Brauer block.

\bibitem[RS]{RagnarssonStancuIdempotents}
K{\'a}ri Ragnarsson and Radu Stancu.
\newblock Saturated fusion systems as idempotents in the double {B}urnside
  ring.
\newblock {\em Geom. Topol.}, 17(2):839--904, 2013.

\bibitem[Th{\'e}]{ThevenazBook}
Jacques Th{\'e}venaz.
\newblock {\em {$G$}-algebras and modular representation theory}.
\newblock Oxford Mathematical Monographs. The Clarendon Press, Oxford
  University Press, New York, 1995.
\newblock Oxford Science Publications.

\end{thebibliography}
\end{document}